\documentclass[pdflatex,sn-mathphys-num]{sn-jnl}

\usepackage{graphicx}%
\usepackage{multirow}%
\usepackage{amsmath,amssymb,amsfonts,bm}%
\usepackage{amsthm}%
\usepackage{amsfonts}
\usepackage{mathrsfs}%
\usepackage[title]{appendix}%
\usepackage{xcolor}%
\usepackage{textcomp}%
\usepackage{manyfoot}%
\usepackage{booktabs}%
\usepackage{algorithm}%
\usepackage{algorithmicx}%
\usepackage{algpseudocode}%
\usepackage{listings}%

 \newcommand{\bs}[1]{\boldsymbol{#1}}

\def \bx{\bs x}
\def \by{\bs y}
\def \bz{\bs z}
 
\def\B{\mathbb{B}}

\def\R{\mathbb{R}}

\def\e{\mathrm{e}}
\def\d{\mathrm{d}}

\def\tX{X}
\def\Is{\mathcal{I}_S^d}
\newcommand \dint {\displaystyle\int}
\def\L2{L^2}

\allowdisplaybreaks[4]

\usepackage{pifont}
\usepackage{subfigure}
\usepackage{comment}

\theoremstyle{thmstyleone}%
\newtheorem{theorem}{Theorem}[section]
\renewcommand{\proofname}{\rm{Proof.}}
\theoremstyle{thmstyletwo}%
\newtheorem{exa}{\bf Example}
\newtheorem{lemma}{\bf Lemma}[section]
\theoremstyle{thmstylethree}%
\newtheorem{remark}{Remark}[section]
\numberwithin{equation}{section} 
\numberwithin{figure}{section} 

\begin{document}

\title[Sparse Monte Carlo method for semi-linear PIDEs]
{Efficient implicit-explicit sparse stochastic method for high dimensional semi-linear nonlocal diffusion equations }


\author*[1]{\fnm{Changtao Sheng}  }\email{ctsheng@sufe.edu.cn}

\author*[2]{\fnm{Bihao Su} }\email{bihaosu@hainanu.edu.cn}

\author[1]{\fnm{Chenglong Xu}  }\email{xu.chenglong@sufe.edu.cn}
 
\affil[1]{\orgdiv{School of Mathematics}, \orgname{Shanghai University of Finance and Economics}, \orgaddress{ \city{Shanghai}, \postcode{200433},  \country{China}}}

\affil[2]{\orgdiv{School of Mathematics and Statistics}, \orgname{Hainan University}, \orgaddress{ \city{Haikou}, \postcode{570100}, \country{China}}}



\abstract{In this paper, we present a sparse grid-based Monte Carlo method for solving high-dimensional semi-linear nonlocal diffusion equations with volume constraints. The nonlocal model is governed by a class of semi-linear partial integro-differential equations (PIDEs), in which the operator captures both local convection-diffusion and nonlocal diffusion effects, as revealed by its limiting behavior with respect to the interaction radius. To overcome the bottleneck of computational complexity caused by the curse of dimensionality and the dense linear systems arising from nonlocal operators, we propose a novel implicit-explicit scheme based on a direct approximation of the nonlinear Feynman-Kac representation. The incorporation of sparse grid interpolation significantly enhances the algorithm’s scalability and enables its application to problems in high dimensions. To further address the challenges posed by hypersingular kernels, we design a sampling strategy tailored to their singular structure, which ensures accurate and stable treatment of the nonlocal operators within the probabilistic framework. Notably, the proposed method inherits unconditional stability from the underlying stochastic representation, without imposing constraints on the temporal and spatial discretization scales. A rigorous error analysis is provided to establish the convergence of the proposed scheme. Extensive numerical experiments, including some non-radial solutions in up to 100 dimensions, are presented to validate the robustness and accuracy of the proposed method.}


\keywords{Brownian motion, Monte Carlo method, error bounds, Sparse grid}


\pacs[MSC Classification]{60G52,  65C05, 65D40, 65L70, 91G60}

\maketitle

\section{Introduction}\label{sec1}

Nonlocal partial differential equations (PDEs), which characterize spatial and temporal interactions beyond local neighborhoods, have been widely applied in fields such as fluid dynamics \cite{del2010}, plasma physics \cite{del2006, Sanchez2008}, turbulence \cite{chen2006, Gunzburger2018}, and materials science \cite{Silling2000,Silling2005}, where peridynamics has been extensively used for crack modeling. Compared to classical local models based on Brownian motion, nonlocal models driven by mixed Lévy processes, which combine diffusive and jump behaviors, offer greater modeling flexibility, feature compact formulations and clearer physical interpretability. Such models have found applications in biological systems \cite{tao2011, Ahmed2007}, image processing \cite{Gilboa2009, Buades2010}, and finance \cite{cont2005integro}. However, even in the linear case, exact solutions for nonlocal models with non-constant kernels are typically unavailable, and when such solutions are known, the corresponding source terms are often difficult to compute explicitly. 

Such nonlocal models are notoriously difficult to solve numerically due to the intrinsic nonlocality of the operators and/or the possible presence of hypersingular integrals. Much effort has been devoted to the development of numerical algorithms for nonlocal models with volume constraints. Some recent studies have focused on the limit of the interaction radius $\delta \to \infty$, which gives rise to PDEs involving the fractional Laplacian. The main deterministic numerical methods employed for these problems include, but are not limited to, finite difference methods~\cite{duo2018finite,hao2019fractional,victor2020simple}, finite element methods~\cite{acosta2017fractional,Ainsworth2017,bonito2019numerical}, and spectral methods~\cite{mao2017hermite,tang2018hermite,sheng2019fast}. Beyond this asymptotic regime, nonlocal models with finite interaction radius $0 \leq \delta < \infty$ have also attracted significant attention, and a variety of promising deterministic numerical methods have been developed in recent years (see, e.g.,~\cite{Tian2013,Leng2021,Aggarwal2024} and references therein). The recent works~\cite{Elia2020,Du2023} and monograph \cite{Du2019} provide comprehensive and up-to-date reviews of the numerical issues associated with various kinds of nonlocal models. However, even in two- and three-dimensional settings, these methods can become highly complex, and their computational cost becomes prohibitive in higher-dimensional cases, mainly due to: (i) the curse of dimensionality, which leads to exponential growth in computational complexity; (ii) the nonlocal nature of the operators, which results in dense system matrices and significantly increases both computational and storage costs; and (iii) the presence of singular kernels, which further complicates numerical approximation.
To address the challenges arising in high-dimensional settings, we also refer the reader to recent advances in deep Neural Network-based methods for solving linear PDEs involving various types of nonlocal operators in multiple dimensions (cf.~\cite{PANG2020109760,Guo2022M,Castro2022,Lu2024}).

In contrast to deterministic methods, particle-based stochastic approaches construct numerical solutions by establishing connections between PDEs and stochastic processes, such as Brownian motion for the classical Laplace operator, $\alpha$-stable Lévy processes for the fractional Laplacian operator, and compound Poisson processes for nonlocal integro-differential operators \cite{DuHuang2014,Metzler2000,Pardoux2005,Robbe2017NM}.  These path-integral-based stochastic methods avoid the construction and solving of the dense stiffness matrices inherent in traditional deterministic schemes for nonlocal problems, and are naturally suited for parallel implementation (cf. \cite{Lovbak2021NM,Shao2020,Ding2023,Lei2025an}). Moreover, when solving high-dimensional problems or those posed on irregular domains, it suffices to simulate stochastic process trajectories within the domain, offering a new perspective for addressing such challenges. In recent years, this class of methods has demonstrated remarkable effectiveness in solving nonlocal problems. For instance, Yang et al.~\cite{yang2023} proposed an insightful probabilistic algorithm based on the Feynman-Kac formula for 3D semi-linear nonlocal diffusion equations with volume constraints, providing a useful approach for tackling such nonlocal problems. Kyprianou et al. \cite{Kyprianou2018} developed a walk-on-spheres Monte Carlo method for solving fractional Poisson equations in 2D. Sheng et al. \cite{SSX2023,SXS2023} further investigated stochastic algorithms for solving fractional PDEs on high-dimensional irregular domains and provided rigorous error analysis. However, to the best of our knowledge, there are no stochastic algorithms with rigorous error estimates available for semi-linear nonlocal diffusion equations with volume constraints in arbitrary high dimensions.

The study of limiting behavior in nonlocal models, in which the nonlocal effects gradually vanish and the corresponding limits naturally bridge fractional and classical local models, has been extensively explored (cf.\cite{Tian2016}). The method proposed in this paper offers a new perspective for understanding such limiting behavior from a stochastic point of view. Specifically, we consider a semilinear nonlocal equation involving advection-diffusion terms, an integral-kernel-based nonlocal operator, and a nonlinear forcing term. We develop a semi-implicit stochastic algorithm based on the nonlinear Feynman-Kac framework, which establishes a probabilistic representation linking PDEs and stochastic processes. Notably, the considered equation (see \eqref{mainprob} in this paper) incorporates both advection-diffusion operators and kernel-based nonlocal operators, naturally leading to a hybrid stochastic process involving both continuous and jump components: the local diffusion corresponds to a continuous process whose trajectories touch the boundary upon exiting the domain, while the nonlocal component corresponds to a jump process that allows trajectories to exit the domain without boundary contact. Therefore, the limiting behavior of these two stochastic processes provides a natural and unified way to capture and connect both the local and nonlocal regimes. 

In this paper, we develop an interpretable and efficient stochastic algorithm for high-dimensional semi-linear nonlocal diffusion equations with volume constraints on bounded domain, which constructs the numerical solution by applying the Feynman-Kac formula over each time subinterval and simulating sample paths of the underlying stochastic process within each subinterval. An implicit-explicit scheme is introduced to treat the nonlinear term, while sparse grids are utilized to reconstruct the numerical solution from the previous time level, thereby substantially reducing the number of degrees of freedom and effectively addressing the computational challenges associated with high dimensionality.  
We highlight below the main contributions of the paper.
\begin{itemize}
\item An efficient stochastic algorithm is proposed for solving semilinear nonlocal diffusion equations with volume constraints, which is developed through a numerical discretization of the nonlinear Feynman-Kac representation and further incorporates sparse grid interpolation and an implicit-explicit time discretization strategy. A distinctive feature of the algorithm is that it relies solely on sample path simulation, thereby entirely bypassing the solution of dense linear systems caused by nonlocal operators, while naturally supporting parallel implementation and exhibiting excellent scalability for high-dimensional nonlocal models. Numerical experiments demonstrate the accuracy and numerical advantages of the proposed method, even in challenging cases with non-radial solutions in up to 100 dimensions.
\item In \cite{yang2023}, due to the adoption of a quadrature-based expectation formulation, the mesh points near the boundary must remain strictly within the domain. This can be achieved by setting the spatial mesh and temporal mesh to satisfy $\text{dist}(x_j, \partial \mathcal{D}) \geq \mathcal{O}\big((\Delta t)^{\frac{1}{2} - \varepsilon}\big),$ $\text{for } x_j \in  \mathcal{D}, $ which ensures that all quadrature points lie within the domain and the numerical scheme remains stable.  
This condition suggests that a proper balance between the spatial and temporal mesh sizes such that  $(\Delta x)^2 \geq \Delta t.$  In contrast, we employ a direct path simulation strategy with sparse grid interpolation applied during the reconstruction step within the implicit-explicit scheme. This formulation effectively circumvents the need for numerical integration of the expectation and inherently eliminates the necessity of imposing additional constraints on the spatial and temporal grids, as further confirmed by the  numerical results.
\item Traditional deterministic methods face significant challenges when dealing with nonlocal operators, particularly those involving hypersingular kernels, whereas stochastic algorithms treat such operators as jump processes and compute them by directly sampling a random variable $\bz_k$ (see \eqref{Xi}).
This implies that the proposed algorithm entails essentially the same computational cost for both local and nonlocal problems, which stands in sharp contrast to the significant disparity observed in traditional deterministic methods.
More importantly, we provide exact sampling formulas for the random variables $\bz_k$ associated with several representative kernel functions, including some hypersingular kernels (see Section~\ref{sec3.2}).

\item We rigorously analyze the error bound under the $L^2$-norm for the proposed implicit-explicit scheme that is totally explicit with respect to the number of simulation paths $M$ and the time step size $\Delta t$. To the best of our knowledge, this is the first work that establishes a rigorous convergence analysis for stochastic methods applied to semi-linear PIDEs in arbitrarily high dimensions.
\end{itemize}

The structure of this article is outlined as follows. Section 2 introduces the governing equations under consideration, derives the Feynman-Kac formula within each time subinterval, and provides a preliminary discussion of relevant background concepts. In Section 3, we present the details of random sampling for several representative kernel functions, the approximation of stochastic trajectories, and briefly review sparse grid interpolation. In Section 4, we establish rigorous error estimates for the proposed algorithm. Section 5 presents extensive numerical experiments to demonstrate its efficiency and accuracy. Finally, concluding remarks are given in Section 6. 

 \section{Preliminaries on main results}
In this section, we state the time-dependent PIDEs of interest along with the necessary assumptions and present its probabilistic representation to establish a rigorous framework for the stochastic algorithm. Lastly, we introduce the time-stepping-based semi-discrete scheme, which lays the foundation for the fully discrete scheme discussed in Section \ref{EIsheme}.
 
 \subsection{Problem setting}
 Let $\Omega \subset \R^{n}$ be a bounded open domain, and let its extended domain  $\Omega_E$ be defined such that $\Omega$ and $\Omega_E$ are disjoint. 
 The mainly concerned in this paper is a general class of semi-linear nonlocal volumetric constrained diffusion equation, which is expressed as
\begin{equation}
\label{mainprob}
\begin{cases}
\partial_{t}u(t,\bx) -\mathcal{L}[u](t,\bx)  =f(t,\bx,u),\;&(t,\bx) \in (0,T]\times \Omega,\\[4pt]
u(t,\bx) = g(t,\bx),\;&(t,\bx) \in (0,T]\times \Omega_{E},\\[4pt]
u(0,\bx) =u_{0}(\bx),\quad &\bx\in  \Omega \cup \Omega_{E},
\end{cases}
\end{equation}
where $f(t, \bx, u)$ (also short for $f(u)$) denotes the forcing term, potentially depending nonlinearly on $u$, and is subject to the volume constraint $u(t, \bx) = g(t, \bx)$ over the nonempty interaction domain $\Omega_{E}$.
In the above, the the non-divergence form integro-differential operator $\mathcal{L}$ in \eqref{mainprob} is defined by
\begin{equation}\label{nonlocaloperator}\begin{split}
\mathcal{L}[u](t,\bx)  =\,&\frac 1 2 \text{Trace}\big( \sigma(t,\bx) \sigma(t,\bx)^{\top} \text{Hess}_{\bx}u(t,\bx) \big) + \big\langle  \mu(t,\bx), \nabla_{\bx}u(t,\bx) \big\rangle_{\R^{d}}\\
& + \dint_{\mathcal{D}}\Big[u(t,\bx+c(t,\bx,\bz)) - u(t,\bx) \Big]\varphi(\bz)\,{\rm d}\bz,
\end{split}\end{equation}
where $\langle \cdot,\cdot \rangle$ denotes the scalar product. The drift term $\mu$ is a known vector-valued function,  the local diffusion coefficient  $\sigma$ is a Lipschitz continuous $d \times d$ matrix-valued function with $\sigma^\top$ denoting its transpose, $c\in \R^d$ denote the jump amplitude vector-valued function, $\nabla_{\bx} u$ and $\text{Hess}_{\bx}u$ denote the gradient and the Hessian of function $u$ with respect to $\bx$, respectively, and $\text{Trace}$ denotes the trace of a matrix. 
The interaction domain $\mathcal{D}\subset \mathbb{R}^d$ for $\bx \in \Omega$ is defined in conjunction with the nonlocal kernel $\varphi(\bz)$, which governs the interactions within the domain.
In this paper, we focus on two classes of kernel functions:  the first one consists of nonnegative and integrable kernels satisfying
\begin{equation}
\label{psi}
\varphi(\bz)\ge 0\; \text{for}\;\bz\in \mathcal{D}\;\;\; \text{and}\;\;\;\psi(\bz) = \frac{\varphi(\bz)}{\lambda},\;\;\;\lambda = \dint_{\mathcal{D}}\varphi(\bz)\,\d\bz<\infty,
\end{equation}
and the second one comprises nonnegative and hypersingular kernels that satisfy
\begin{equation}
\label{lambda}
\varphi(\bz)\ge 0\; \text{for}\;\bz\in \mathcal{D}\;\;\; \text{and}\;\;\;\psi(\bz) = \frac{\varphi(\bz)|\bz|^2}{\lambda},\;\;\;\lambda =\int_{\mathcal{D}} \varphi(\bz)|\bz|^{2}\,\d\bz<\infty,
\end{equation}
where $\psi(\cdot)$ can be viewed as a probability density function. 
Now we are ready to show that the relationship between the extended domain $\Omega_{E}$ and the interaction domain $\mathcal{D}$ is as follows:
\begin{equation}
\Omega_{E} = \big\{\by\in \Omega^{c}: \by=\bx+c(t,\bx,\bz), \;\; \text{for} \;\; \bz\in\mathcal{D}, \;\bx\in \Omega,\; t\in[0,T] \big\}.
\end{equation}
When $c(t,\bx,\bz)=\bz$ and $\mathcal{D}=\B^d_\delta$, which is regarded as a specific case in this study, the extended domain reduces to $\Omega_E=\left\{\bs x \in \mathbb{R}^{d} \backslash \Omega\,:\, \text {dist}(\bs x, \partial \Omega) \leq \delta\right\}$ with the radius $\delta > 0$, a scenario of particular interest to us that will be discussed in detail in Section \ref{sec3.2}.
The nonlocal component in $\mathcal{L}$ corresponds to the compound Poisson process. In the particular case where the kernel reduces to a constant, the operator $\mathcal{L}$ is governed by a process comprising a mixture of Brownian motion and a compound Poisson process. Conversely, if the kernel is given by $\varphi(\bz) = c_{d,\delta}^{\alpha} |\bz|^{-d-\alpha}$ and both the drift and diffusion terms in $\mathcal{L}$ vanish, the problem \eqref{mainprob} simplifies to a fractional diffusion equation driven by an $\alpha$-stable Lévy process. On the other hand, when the drift and diffusion terms are nonzero, the operator $\mathcal{L}$ is driven by a process that combines Brownian motion and an $\alpha$-stable Lévy process. The well-posedness of the problem described in \eqref{mainprob} has been established in \cite{DuGunzburger2012} under standard assumptions.

\subsection{The Feynman-Kac formula}
This subsection presents the Feynman-Kac formula associated with semi-linear PIDEs on bounded domains \cite{Zhang2012}, which serves as the theoretical foundation for the development of the stochastic algorithms in this paper.

Let $\{W_{t}\}_{t\in(0,T]}$ denote a $d$-dimensional standard Brownian motion, and let $\{X_{t}\}_{t\in(0,T]}$ be an adapted  $d$-dimensional stochastic process with continuous sample paths, which satisfies the following stochastic equation
\begin{equation}
\label{Xt}
X_{t} = \bx + \int_{0}^{t} \mu(s, X_{s})\,\d s + \int_{0}^{t}\sigma(s,X_{s})\,\d W_{s} + \sum_{k=1}^{N_{t}}c(t,\bx,\bz_{k}),
\end{equation}
where the coefficients $\mu,\sigma,c$ are defined in \eqref{nonlocaloperator} and $N_{t}$ is a Poisson process with intensity determined by the corresponding Poisson probability distribution
\begin{equation}
\label{PNk}
\mathbb{P}(N_{t} = k) = \frac{(\lambda t)^{k}}{k!}e^{-\lambda t},
\end{equation}
with $\lambda$ defined in \eqref{psi} (resp. \eqref{lambda}). 
Thus, the Feynman-Kac formula for \eqref{mainprob} admits a unique solution with a stochastic representation, and for every $\bx\in\R^{n}$, it holds that
\begin{equation}
\begin{split}
\label{FK}
u(t,\bx) = \mathbb{E}_{X_0=\bx}\Big[ u_{0}(X_{t}) \mathbb{I}_{\tau_{\Omega}>t}  + g(X_{\tau_{\Omega}},\tau_{\Omega}) \mathbb{I}_{\tau_{\Omega}\leq t} \Big]  + \mathbb{E}_{X_0=\bx}\Big[  \int_{0}^{t \land \tau_{\Omega}} f(u(s,X_{s}))\,\d s  \Big],
\end{split}\end{equation}
where $X_{t}$ represents a stochastic process satisfying \eqref{Xt} with $X_{0}=\bx$, $\tau_{\Omega} = \inf \{ s| X_{s} \notin \Omega \}$ represents the first exit time from $\Omega$, and $\mathbb{I}$ is the indicative function.

We observe that the formula \eqref{FK} establishes a direct connection between the solution and stochastic processes, which implies that the first step in developing efficient stochastic algorithms is to solve \eqref{Xt} or, equivalently, simulate the paths of the stochastic process. Fortunately, if \( X_0 \) is assumed to follow a known single-point distribution, \( X_0 = \xi \in \mathbb{R}^d \), the entire process \( X_t \) can be fully characterized and simulated. However, the first challenge we encounter is that the Feynman-Kac formula \eqref{FK} is nonlinear, which may require an iterative procedure to obtain a numerical solution. This, however, can result in a significant computational burden, making such an approach less favorable. To address this, we will introduce a time-stepping method in the forthcoming subsection to overcome this issue.

\subsection{The time stepping scheme}\label{sec2.3}

To establish the time-stepping scheme, we first construct a uniform grid in the time domain $[0,T]$ to simplify the presentation, where the time interval is evenly divided into $N$ parts as follows
\begin{equation}
\label{T}
\mathcal{T} := \{ 0 = t_{0}<t_{1}<t_{2}<\cdots<t_{N-1}<t_{N} = T \},
\end{equation}
and the time step size $t_{i}-t_{i-1}\equiv\Delta t = T/N$, $\mathcal{T}_i=(t_{i-1},t_i],$ $i=1,2,\cdots,N$. Then, it is necessary to further investigate the position of the stochastic process $X_s$ at the time grid points and determine whether it exits the domain $\Omega$. To this end, we define the exit time of the process $X_t$ from the domain $\Omega$ as

\begin{equation}
\tau_{i} := \inf \big\{  s : X_{s}\notin \Omega ,t_{i-1}<s\leq t_{i} \big\},\,\,\, 1\leq i\leq N.
\end{equation}
Within each subinterval $[t_{i-1}, t_{i}]$, it follows from \eqref{Xt} that the stochastic process initiates at $X_{t_{i-1}}$ and transitions to $X_{t_i}$ as
\begin{equation}
\label{Xi}
X_{t_{i}} = X_{t_{i-1}}  + \int_{t_{i-1}}^{t_{i}\land \tau_{i}} \mu(s, X_{s})\,\d s + \int_{t_{i-1}}^{t_{i}\land \tau_{i}} \sigma(s,X_{s})\,\d W_{s} + \sum_{k=1}^{N_{t_{i}\land \tau_{i} -t_{i-1}}}c(t,\bx,\bz_{k}),
\end{equation}
where $s\in(t_{i-1},t_{i}\land \tau_{i}]$ and the coefficients $\mu,\sigma,c$ are defined in \eqref{nonlocaloperator}.  Then, the most natural approach is to discretize the above equation directly using the Euler-Maruyama approximation scheme  
\begin{equation}
\label{Xidiscrete}\begin{split}
\tX_{t_{i}} \approx\;& \tX_{t_{i-1}}  + \mu(t_{i-1}, \tX_{t_{i-1}} ) (t_{i}\land \tau_{i} -t_{i-1}) \\&+ \sigma(t_{i-1},\tX_{t_{i-1}} ) (W_{\tau_{i} \land t_{i}} -W_{t_{i-1}}) + \sum_{k=1}^{N_{\Delta t}}c(t_{i-1},\tX_{t_{i-1}},\bz_{k}).
\end{split}\end{equation}
According to \cite{Higham2005}, the approximation symbol in the above equation can be replaced by an equality in the sense of expectation as $\Delta t \to 0$.

Next, we proceed with the derivation of the time-stepping scheme based on the Monte Carlo method. To this end, we denote the solution at the $i$-th subinterval as follows:
$$u^{i}(t,\bx) = u(t,\bx),\;\;\; \forall t\in\mathcal{T}_i=(t_{i-1},t_i],\;\;\;1\leq i\leq N.$$
Utilizing a time-stepping method, we aim to iteratively solve the underlying problem within each interval $\mathcal{T}_{i}$. By breaking down the problem into smaller, manageable segments, we can effectively leverage the Monte Carlo method to approximate the solution. This methodology not only enhances computational efficiency but also allows us to capture the stochastic nature of the processes more accurately. Thus, we will now outline the specific problem we intend to address on the interval $\mathcal{T}_{i}$:
\begin{equation}
\label{pi}
\begin{cases}
\partial_{t}u(t,\bx) -\mathcal{L}[u](t,\bx)  =f(u),\;&(t,\bx) \in (t_{i-1},t_{i}]\times \Omega,\\[4pt]
u(t,\bx) = g(t,\bx),\;&(t,\bx) \in (t_{i-1},t_{i}]\times \Omega_{E},\\[4pt]
u(t_{i-1},\bx) =u^{i-1}(t_{i-1},\bx),\quad &\bx\in  \Omega \cup \Omega_{E}.
\end{cases}
\end{equation}
Thus, according to \cite{barles1997backward,pardoux1990adapted}, the Feynman-Kac representation of \eqref{pi} within the interval $\mathcal{T}_i$ can be derived as follows
\begin{equation}
\begin{split}
\label{FKi}
u(t_{i},\bx) = \;&\mathbb{E}_{X_{t_{i-1}}=\bx}\Big[ u^{i-1}(t_{i-1},X_{t_{i}}) \mathbb{I}_{\tau_{i}>t_{i}}  + g(\tau_{i},X_{\tau_{i}}) \mathbb{I}_{\tau_{i}\leq t_{i}} \Big] \\& + \mathbb{E}_{X_{t_{i-1}}=\bx}\Big[  \int_{t_{i-1}}^{t_{i} \land \tau_{i}} f(u(s,X_{s}))\,\d s  \Big],
\end{split}\end{equation}

\begin{remark} {\em 
Unlike existing studies that employ numerical discretization of backward stochastic differential equations {\rm\cite{yang2023}} or forward-backward stochastic differential equations {\rm\cite{Lu2024}} to solve PIDEs, our approach adopts a more direct and conceptually streamlined methodology. Specifically, we construct the stochastic process paths $X_t$ by simulating them directly via {\rm\eqref{Xidiscrete}}. This approach remains to determine $\bz_k$ to generate the complete path $X_t$, as detailed in forthcoming {\rm Section \ref{sec3.2}}.
}
\end{remark}
 
 \section{Fully discrete Monte Carlo scheme based on sparse grid.}\label{EIsheme}
This section focuses on the detailed numerical implementation of the time-stepping scheme introduced in the previous section \ref{sec2.3}. In particular, we address two key challenges that remained unresolved: (i).\,the treatment of the integral terms, specifically the random sampling of $\bz$ in \eqref{Xidiscrete}; and (ii).\,the handling of the nonlinear terms, achieved through the implicit-explicit scheme combined with interpolation technique. To efficiently manage high-dimensional problems, we adopt sparse grid interpolation. The section concludes with a comprehensive summary of our algorithm.

\subsection{Treatment of the kernel functions}\label{sec3.2}
The expression of kernel functions plays a critical role in determining the characteristics of the associated jump processes. If the kernel function $\varphi(z)$ is bounded and possesses compact support, numerical quadrature methods such as the Gauss-Legendre or Newton-Cotes rules can be effectively employed. Conversely, for kernel functions singular at the origin, such as $\varphi(z) =|z|^{-\gamma}$ with $0 < \gamma< 1$, the Gauss-Jacobi quadrature rule is a more appropriate choice (cf.\,\cite{yang2023}). Kaliuzhnyi-Verbovetskyi et al. \cite{Kaliuzhny2022} also investigated singular kernel functions of the form $W(\bx,\by) = |\bx-\by|^{-\gamma}$, where $0 < \gamma < d/2$.
This subsection focuses on the case where the interaction domain is $\mathcal{D}=\mathbb{B}^d_\delta$, which has garnered significant attention in the context of peridynamics. Noted that for more general interaction domains $\mathcal{D}=\Omega_\delta = \{\bs x \in \mathbb{R}^{d} \backslash \Omega\,:\, \text {dist}(\bs x, \partial \Omega) \leq \delta\}$ with kernel function satisfies $\varphi(\bx,\by)={\varphi}(|\bx-\by|)$, analogous methods can be applied by introducing a coordinate transformation to map the domain into a spherical domain $\mathbb{B}^d_\delta$.
The subsequent discussion highlights the treatment of several representative classes of kernel functions, which include:
\begin{itemize}
\item[\ding{172}] hypersingular kernel functions: $C_H|\bz|^{-d-\alpha}$, with $C_H$ denotes a normalization constant;
\item[\ding{173}] tempered hypersingular kernel functions: $C_T {\rm e}^{-\beta|\bz|}|\bz|^{-d-\alpha}$, with $C_T$ denotes a normalization constant;
\item[\ding{174}] Gaussian kernel functions: $C_G {\rm e}^{-\frac{|\bz|^2}{2\sigma^2}}$, with $C_G$  denotes a normalization constant;
\end{itemize}

A notable property of the nonlocal operator, when it is reduced to the form $\mathcal{L}_\delta u(\bx)= \int_{\mathbb{B}^d_\delta}(u(\bx+\bz)-u(\bx))\varphi(|\bz|)\d \bz$ and $\varphi(|\bz|)$ is associated with hypersingular kernel functions \ding{172}, is its convergence to $-\Delta$ as $\delta \to 0$ and to $(-\Delta)^{\frac{\alpha}{2}}$ as $\delta \to \infty$.  This behavior underscores its role as a transitional model seamlessly connecting integer-order models and fractional-order models, see \cite{Tian2016} for a detailed discussion. 
\begin{lemma}\label{lemker1}{\bf (Hypersingular kernel).} {\em
Let $\alpha\in(0,2)$, and defined the hypersingular kernel functions as
\begin{equation}\label{kernelsingular}
\varphi_H(\bz) =C_{H} |\bz|^{-d-\alpha}, \;\;\;C_H:=C_H(d,\delta,\alpha) = (2-\alpha)\delta^{\alpha-2}/\omega_d,\;\;\; \bz\in \mathbb{B}^d_\delta,
\end{equation}
where $\omega_d=2\pi^{d/2}/\Gamma(d/2)$ denotes the surface area of the unit sphere $\mathbb{S}^{d-1}=\{\bx\in\mathbb{R}^d:|\bx|=1\}$, and the kernel function satisfies the bounded second order moment condition
\begin{equation}\label{secondmoment}
\int_{\B_{\delta}^d} \varphi_H(|\bz|)|\bz|^{2}\,\d\bz = 1.
\end{equation}
Denote the unit vector along any nonzero vector $\bz$ by $\hat \bz =\bz/r$ with $r=|\bz| := \sqrt{ \langle\bz, \bz\rangle}$, then the stochastic sampling of $\bz$ in the nonlocal term involving the hypersingular kernel function \eqref{kernelsingular} can be computed as $\bz=r\hat{\bz}$, where 
\begin{equation}
r = \delta\,\xi^{1/(2-\alpha)},\;\;  \hat{\bz} = \frac{(z_1, z_2, \dots, z_d)}{\sqrt{z_1^2 + z_2^2 + \dots + z_d^2}},\;\;\xi\sim U(0,1),\;\;\bz\sim U(0,1)^d.
\end{equation}}
 \end{lemma}\vspace{-20pt}
\begin{proof}
We commence by  verifying whether the kernel function $\varphi_H(|\bz|)|\bz|^2$ fulfills the requirements of a probability density function, namely that $\int_{\B_\delta^d} \varphi_H(|\bz|)|\bz|^2 \d\bz = 1$.
Using the spherical volume element $\d\bz =r^{d-1} \d r\, \d\sigma (\hat\bx)= \omega_d r^{d-1} \d r$ and noting that $C_H= (2-\alpha)\delta^{\alpha-2}/\omega_d$, we find from \eqref{kernelsingular}  that
\begin{equation}
\begin{split}
\int_{\B_\delta^d} \varphi_H(|\bz|)|\bz|^2 \,\d\bz &= 
C_H\omega_d \int_0^\delta r^{-\alpha+1} \,\d r=C_H\frac{\delta^{2-\alpha} \omega_d}{2-\alpha}=1,
\end{split}
\end{equation}
which confirms that $\varphi_H(|\bz|) |\bz|^2$ constitutes a valid probability density function.
In view of the fact that $\varphi_H(|\bz|)|\bz|^2\d\bz =\omega_d C_H  s^{1-\alpha} \d s=\omega_d\, C_H  s^{1-\alpha}\d s$, then the cumulative distribution function $F(r)$ can be computed by
\begin{equation}
\begin{split}
F(r) =& \int_{\B_r^d} \varphi_H(|\bz|) |\bz|^{2} \,\d \bz= \omega_d C_H  \int_0^r s^{1-\alpha}\,\d s =  \frac{2-\alpha}{\delta^{2-\alpha}}\int_0^rs^{1-\alpha}\,\d s
= \left(\frac{r}{\delta} \right)^{2-\alpha}.
\end{split}
\end{equation}
By applying the inverse transform sampling method and setting $F(r) = \xi$ with $\xi \sim U(0,1)$, we solve for $r$:
\begin{equation}
r = F^{-1}(\xi) = \delta \,\xi^{1/(2-\alpha)}.
\end{equation}
What remains is to generate a uniform distribution on the $d$-dimensional unit sphere $\hat{\bz} = \bz / |\bz|$. To achieve this, we first generate $d$ independent standard normal random variables $\bz\sim U(0,1)^d$ (indicating $z_i$ are i.i.d. on $[0,1]$) and then compute $\hat{\bz}$ by normalizing $\bz$ as
$$  \hat{\bz} = \frac{(z_1, z_2, \dots, z_d)}{\sqrt{z_1^2 + z_2^2 + \dots + z_d^2}}.$$ Finally, by combining the radial component $r$ with the directional vector $\hat{\bz}$, the sampled point is expressed as \(\bz = r \,\hat{\bz}\), yielding the desired result. 
   \end{proof}

   \begin{lemma} \label{lemker2}{\bf (Tempered hypersingular kernel).} {\em Let $\alpha\in(0,2)$, $\beta>0$, and defined the tempered hypersingular kernel functions as
\begin{equation}\label{temperedkernel}
\varphi_{T}(\bz) = C_{T}\e^{-\beta|\bz|} |\bz|^{-d-\alpha},\;\;\;C_T:=C_T(d,\delta,\alpha,\beta)=\frac{\beta^{2-\alpha}}{\omega_d \gamma(2-\alpha, \beta \delta)},\;\;\; \bz\in \mathbb{B}^d_\delta,
\end{equation}
where $\gamma(s, x)$ represents the lower incomplete gamma function and $C_T$ denotes a normalization constant chosen to satisfy the bounded second order moment condition \eqref{secondmoment}. 
 Then, the random variable $\bz$ in the nonlocal term involving the hypersingular kernel function \eqref{temperedkernel}  can be computed as $\bz = r\hat{\bz}$, where
 \begin{equation}
r =   \gamma^{-1}\big(2-\alpha, \xi \, \gamma(2-\alpha, \beta \delta )\big), \;\;\; \hat{\bz} = \frac{(z_1, z_2, \dots, z_d)}{\sqrt{z_1^2 + z_2^2 + \dots + z_d^2}}, \;\;\; \xi\sim U(0,1),\;\;\;\bz\sim U(0,1)^d, 
\end{equation}
with $\gamma^{-1}(s,y)$ denotes the inverse of lower incomplete gamma function.}
 \end{lemma}\vspace{-20pt}
\begin{proof}
We first verify that the normalization constant satisfies the second-order moment condition.
Using the spherical volume element $\d\bz =s^{d-1} \d s\, \d\sigma (\hat\bx)= \omega_d s^{d-1} \d s$, where $\omega_d = \frac{2\pi^{d/2}}{\Gamma(d/2)}$ we obtain from \eqref{temperedkernel} that
\begin{equation}
\begin{split}
\int_{\B_{\delta}^{d}} \varphi_{T}(\bz)|\bz|^2 \, \d\bz  
&=C_T \omega_d \int_0^\delta e^{-\beta s} s^{1-\alpha} \, \d s=C_T \omega_d\int_{0}^{\beta\delta}e^{-\zeta} \left( \frac{\zeta}{\beta} \right)^{1-\alpha}\frac{\d \zeta}{\beta}
\\& =C_T \omega_d \beta^{\alpha-2}\int_{0}^{\beta\delta}e^{-\zeta} \zeta^{1-\alpha}\,\d \zeta =C_T \omega_d \beta^{\alpha-2}\gamma(2-\alpha,\beta \delta)=1.
\end{split}
\end{equation}
Similarly, the cumulative distribution function can be expressed as
\begin{equation}
\begin{split}
F(r) = \int_{\B_{r}^{d}} \varphi_{T}(\bz) |\bz|^2\d\bz = C_T \omega_d \beta^{\alpha-2} \gamma(2-\alpha, \beta r) = \frac{\gamma(2-\alpha, \beta r)}{\gamma(2-\alpha, \beta \delta)} := \xi \in[0,1],
\end{split}
\end{equation}
which implies
\begin{equation}
\gamma(2-\alpha, \beta r) = \xi\gamma(2-\alpha, \beta \delta) \Longrightarrow r =  \gamma^{-1}\big(2-\alpha, \xi \gamma(2-\alpha, \beta \delta )\big).
\end{equation}
This ends the proof.
   \end{proof}
   
 \begin{lemma}\label{lemker3} {\bf (Gaussian kernel).} {\em
Denote the Gaussian kernel function by  
\begin{equation}\label{guasskernel}
\varphi_G(|\bz|) = C_{G} \e^{-\frac{|\bz|^2}{2\sigma^2}},\;\;\;C_G:=C_G(d, \delta,\sigma)= \frac{2^{1-d/2}}{\omega_d \sigma^d \gamma\left(\frac{d}{2}, \frac{\delta^2}{2\sigma^2}\right)},\;\;\; \bz\in \mathbb{B}^d_\delta,
\end{equation}
where $C_G$ is the normalization constant ensuring that $\int_{\B_\delta^n} \varphi_G(|\bz|) \d\bz = 1$, and $\sigma > 0$ is the standard deviation of the kernel. Then, the random variable $\bz$ in the nonlocal term involving the smooth kernel function \eqref{guasskernel}  can be computed as $\bz = r\hat{\bz}$, where
\begin{equation}
r=\sigma\sqrt{2\gamma^{-1}\Big(\frac{d}{2}, \xi\, \gamma\Big(\frac{d}{2}, \frac{\delta^2}{2\sigma^2}\Big)\Big)},\quad \hat{\bz} = \frac{(z_1, z_2, \dots, z_n)}{\sqrt{z_1^2 + z_2^2 + \dots + z_n^2}}\quad \xi\sim U(0,1),\;\;\bz\sim U(0,1)^d.
\end{equation}
}
\end{lemma}\vspace{-20pt}
\begin{proof}
We begin by verifying the probability density function, which can be derived from \eqref{guasskernel} as follows
\begin{equation*}\begin{split}
\int_{\B_\delta^d} \varphi_G(|\bz|)\,\d\bz &=C_G\,\omega_d \int_0^\delta  \e^{-\frac{r^2}{2\sigma^2}} r^{d-1}\,\d r
\\&=2^{d/2-1}\,C_G\, \omega_d \,\sigma^d \int_0^{\frac{\delta^2}{2\sigma^2}} \e^{-\eta}\eta^{(d/2)-1}  \, \d\eta=2^{d/2-1}\,C_G \,\omega_d \,\sigma^d \gamma\Big(\frac{d}{2},\frac{\delta^2}{2\sigma^2}\Big)=1.
\end{split}\end{equation*}
Similarly, we can express the corresponding cumulative distribution function as:  
\begin{equation*}
\begin{split}
F(r) &= C_G\,\omega_d\int_{0}^{r}  \e^{-\frac{\zeta^2}{2\sigma^2}}\zeta^{d-1}\,\d \zeta=2^{d/2-1}\,C_G \,\omega_d \,\sigma^d \gamma\Big(\frac{d}{2},\frac{r^2}{2\sigma^2}\Big)=\frac{\gamma\left(\frac{d}{2}, \frac{r^2}{2\sigma^2}\right)}{\gamma\left(\frac{d}{2}, \frac{\delta^2}{2\sigma^2}\right)}:=\xi\sim U[0, 1],
\end{split}
\end{equation*}
which directly implies the desired result. This complete the proof.
\end{proof}
Note that the lower incomplete gamma function $\gamma(s, x)$ and its inverse $\gamma^{-1}(s, y)$, which are required for the sampling process, can be directly computed using MATLAB's built-in functions \texttt{gammainc.m} and \texttt{gammaincinv.m}.
\begin{remark} {\em 
Unlike the Gauss quadrature-based approach proposed in \cite{yang2023}, which is efficient in low-dimensional settings but becomes computationally expensive in high dimensions or non-spherical domains, our method directly simulates the underlying stochastic process. This strategy simplifies implementation and improves the handling of nonlocal terms. When the probability density associated with the hypersingular kernel admits a closed-form inverse cumulative distribution function (CDF), the jump length can be computed exactly. Otherwise, standard sampling techniques can approximate the jumps with high accuracy, provided that the sample size is sufficiently large. This approach offers a favorable balance between efficiency and accuracy, particularly in high-dimensional applications.}
\end{remark}
\begin{remark}{\em
If the integration domain is irregular, i.e., replacing the domain $ \B^d_{\delta} $ with a more general domain $\mathcal{D} $, Monte Carlo methods can be employed for random sampling. To accommodate irregular integration domains, several techniques are available:
\begin{itemize}
\item \text{Rejection Sampling}: This technique samples from a simpler distribution $g(x)$ and accepts or rejects samples according to the probability ratio $\frac{f(x')}{M g(x')}$, thereby guaranteeing alignment with the target domain. 
\item \text{Importance Sampling}: This method transforms the integration domain into a simpler one through an appropriate choice of probability distribution, which simplifies the integral evaluation via transformed variables. Importance sampling proves particularly effective for irregular, high-dimensional domains (see, e.g., {\rm\cite{Eriksson2011NM}}). 
\item \text{Markov Chain Monte Carlo (MCMC)}: For more complex, high-dimensional problems, MC -MC methods construct a Markov chain whose stationary distribution approximates the target distribution. This approach enables efficient sampling in intricate and high-dimensional domains. 
\end{itemize}
The choice of method depends on the dimensionality of the problem and the shape of the domain. For low-dimensional and relatively simple regions, rejection sampling is often sufficient. However, for high-dimensional or complex domains, a combination of importance sampling and MCMC techniques may be necessary. These methods are particularly well-suited for sampling the jump lengths of a stochastic process when jumps occur, providing robust solutions for irregular integration problems.}
\end{remark}

\subsection{The approximation of the trajectory $\{\tX_{t_i}\}$}\label{sec3.3}
Next, we turn to the concrete implementation of the stochastic paths $\tX_{t_i}$ in \eqref{Xidiscrete}. 
Since the Feynman-Kac formula \eqref{FKi} incorporates the first exit time $\tau_i$ of the stochastic process $X_t$ from the domain $\Omega$, it is essential to determine whether $X_t$ exits the region $\Omega$ within the time interval $[t_{i-1}, t_i]$. If $X_t$ does exit the domain during this interval, further analysis is required to establish whether the exit is attributable to a jump event. Moreover, the probability of the Poisson process $N_t$ experiencing $k$ jumps is characterized by \eqref{PNk}, as follows:
$$\sum_{k=2}^{\infty} \mathbb{P}(N_{\Delta t} = k) = \mathcal{O}((\Delta t)^2).$$

Therefore, in our calculations, we only need to consider the cases $k = 0$ and $k = 1$. To this end, we define the process $X_t$ to start at $(t_{i-1}, \tX_{t_{i-1}})$ and terminate at $(t_i \land \tau_i, \tX_{t_i \land \tau_i})$.
For the stochastic algorithm, it remains to determine whether the points $X_{t_{i}}$ are inside the domain. More precisely, if $N_{\Delta t} = 0$,  it indicates that there are no jumps in the  process during the time interval, $X_{t}$ is solely driven by Brownian motion. Conversely, when $ N_{\Delta t} = 1 $, it signifies that the process $ X_t $ experiences a single jump during the designated time interval. In this case, $ X_t $ is influenced by two distinct components: the continuous evolution driven by Brownian motion and the discrete impact of a Poisson jump. To simplify the algorithm, we approximate the first arrival time \(\tau_i\) using \(t_i\) in our calculations, due to the complexity of directly simulating $\tau_i$. Furthermore, since the stochastic process exits the domain $\Omega$ at time $\tau_i$ and ceases to move thereafter, we have $\tX_{t_i} = \tX_{\tau_i}$. Consequently, the process \(\tX_{t_i}\) can be formulated as
\begin{equation}
\tX_{t_{i}} = \tX_{t_{i-1}}  + \mu(t_{i-1}, \tX_{t_{i-1}}) \Delta t + \sigma(t_{i-1},\tX_{t_{i-1}}) W_{\Delta t} + \sum_{k=1}^{N_{\Delta t} \land 1}c(t_{i-1},\tX_{t_{i-1}},\bz_{k}),
\end{equation}
where $\bz_k$ is obtained using the method described in Section \ref{sec3.2}.
 While this approach simplifies the computational process, it also introduces some inherent errors, which we will analyze in the next section.


\subsection{Sparse grid interpolation}\label{sec3.1}
To advance the time-stepping scheme for nonlinear equations, we note from \eqref{pi} that the numerical solution at $t_i$ depends on that of the preceding step $t_{i-1}$, regardless of whether an explicit or fully implicit scheme is used.
Owing to the randomness of the sample path locations for stochastic methods, the solution must be reconstructed at arbitrary spatial position from the known discrete data, necessitating the construction of an appropriate interpolation function.
This reconstruction step is critical to the overall efficiency of the proposed algorithm.
To address the computational challenges posed by high-dimensional problems, we introduce sparse grid methods in this section.

Sparse grid methods are designed to mitigate the curse of dimensionality commonly encountered in traditional grid-based approaches (cf.\cite{Bungartz2004NChief,Griebel1999NChief}). The concept behind these methods can be traced back to Smolyak \cite{Smolyak1963}, who initially introduced it for numerical integration. These methods employ a hierarchical basis constructed from nested Chebyshev-Gauss-Lobatto quadrature, as reviewed by Yserentant \cite{Yserentant1986, Yserentant1992}. This approach provides two key advantages: (a) it is spectrally accurate with nested points; (b) FFT can be used for the transform. The sparse grid strategy significantly reduces the number of grid points by selecting only a subset of the full tensor-product grid, while preserving approximation accuracy. In addition, various related methods and applications have been developed, such as those that incorporate structural sparsity, combinatorial selection strategies, and hybrid discretization techniques (see, e.g., \cite{Temlyakov1986,SY2010,SW2010,Harbrecht2008NM,Garcke2006,zhang2013,Zhang2016NM,Adcock2022NM}).

For $d=1$, if $ v $ is a smooth function with a single scalar input $ x $, its interpolation at level $l$ can be written as
\begin{equation}
\label{interp1}
\mathcal{I}^{(l)}[v](x) = \sum_{j=0}^{m_{l}} v(x^{(l)}_j) a^{(l)}_j(x),
\end{equation}
where $x^{(l)}_j $ are the interpolation nodes and $ a^{(l)}_j(x) $ are the corresponding basis functions. The integer $l$, referred to as the interpolation level, indexes a particular one-dimensional set of collocation points $ \{x^{(l)}_{i}\} $, for $ j = 1, \ldots, m_l $. As the level increases, more interpolation points are incorporated into the construction of \eqref{interp1}. Typically, these points correspond to the abscissas of a selected quadrature rule, such as Chebyshev-Gauss-Lobatto points .

 In the multi-dimensional scenario $(d> 1)$, the one-dimensional interpolation and quadrature rules can extend  into $d$ dimensions using tensor products. ns. For this purpose, we introduce a $d$-dimensional multi-index $ \bm{l} := (l_1, l_2, \ldots, l_d) $, where each component $l_i \geq 1$ denotes the level of the one-dimensional quadrature rule associated with the $i$-th input variable. As discussed earlier, tensor-product constructions of one-dimensional quadrature rules are utilized to address multi-dimensional problems. In this context, the tensor-product interpolant corresponding to a multi-index  $\bm{l}$ is given by 
\begin{equation}
\label{interpn}
\mathcal{I}^{(\bm{l})}[v](\bx) = (\mathcal{I}^{l_{1}} \otimes \cdots \otimes \mathcal{I}^{l_{d}})[v](\bx) = \sum_{j_1=1}^{m_{l_{1}}} \cdots \sum_{j_n=1}^{m_{l_{d}}} v(x_{j_1}^{l_{1}}, \ldots, x_{j_d}^{l_{d}}) \cdot (a_{j_1}^{l_{d}} \otimes \cdots \otimes a_{j_d}^{l_{d}}),
\end{equation}
where the number of grid points at level $l$ is defined as
\begin{equation*}
x_j^{(l)} =
\begin{cases}
0, & \text{for } j = 1, \text{ if } m_l = 1, \\[2pt]
-\cos\left(\frac{\pi(j - 1)}{m_l - 1}\right), & \text{for } j = 1, \ldots, m_l, \text{ if } m_l > 1,
\end{cases}\;\;\;\;m_l =
\begin{cases}
1, & \text{if } l = 1, \\
2^{l-1} + 1, & \text{if } l > 1.
\end{cases}
\end{equation*}
For a full tensor-product grid with identical resolution in each dimension (i.e., $m_{l_1} = m_{l_2} = \cdots = m_{l_d}$), the total number of required function evaluations equals $\prod_{i=1}^d m_{l_i}$. This number increases exponentially with dimension $d$, which renders such approaches computationally infeasible in high-dimensional settings. Each additional dimension introduces a multiplicative growth in cost, which leads to the well-known curse of dimensionality. To mitigate this explosion in cost, the Smolyak interpolation method \cite{Barthelmann2000} constructs the interpolant as a linear combination of tensor-product interpolants defined on a hierarchy of grids with varying resolution in different directions. This approach, also known as sparse-grid interpolation, dramatically reduces the number of nodes by leveraging coarse grids in some dimensions while using finer grids only where needed. The resulting sparse-grid interpolant at overall level $m$ can be written as

\begin{equation}
    \label{sgqr}
\Is[v](\bx) = \sum_{|\bm{l}| = m -d + 1}^{m} (-1)^{m - |\bm{l}|} \binom{d-1}{|\bm{l}|}  (\mathcal{I}^{l_{1}} \otimes \cdots \otimes \mathcal{I}^{l_{d}})[v](\bx) ,
\end{equation}
where $m \geq d $, $ |\bm{l}| = l_1 + l_2 + \cdots + l_d $. 

\begin{figure}[!h]\hspace{0pt}
\subfigure{
\begin{minipage}[t]{0.4\textwidth}
\centering
\rotatebox[origin=cc]{-0}{\includegraphics[width=1.1\textwidth,height=1\textwidth]{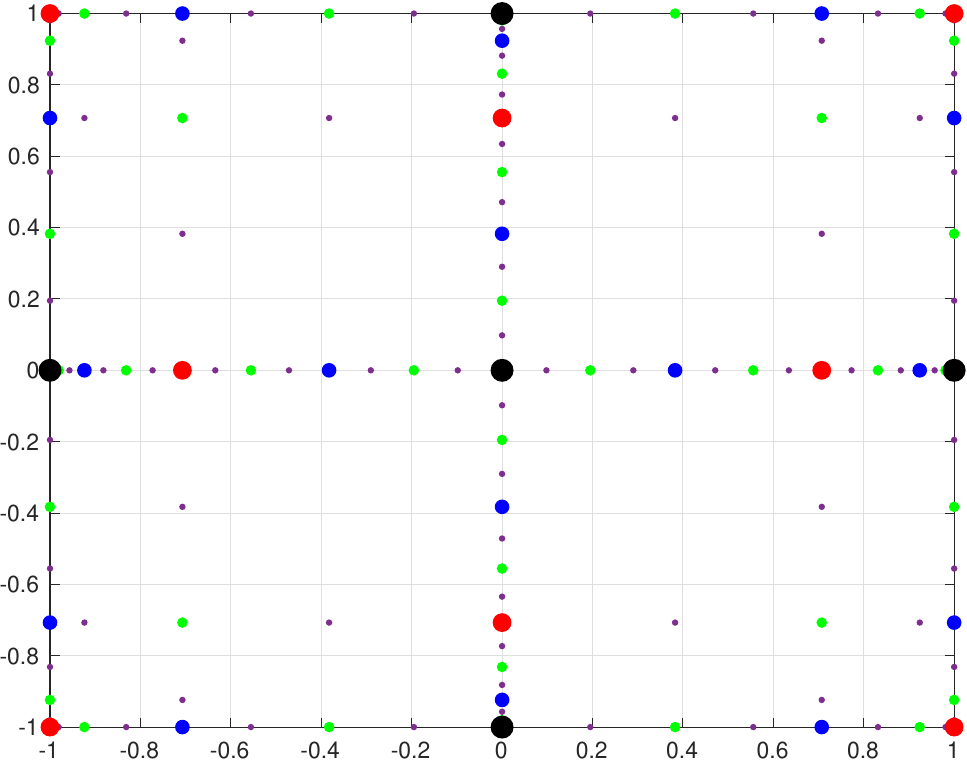}}
\end{minipage}}\hspace{35pt}
\subfigure{
\begin{minipage}[t]{0.4\textwidth}
\centering
\rotatebox[origin=cc]{-0}{\includegraphics[width=1.1\textwidth,height=1\textwidth]{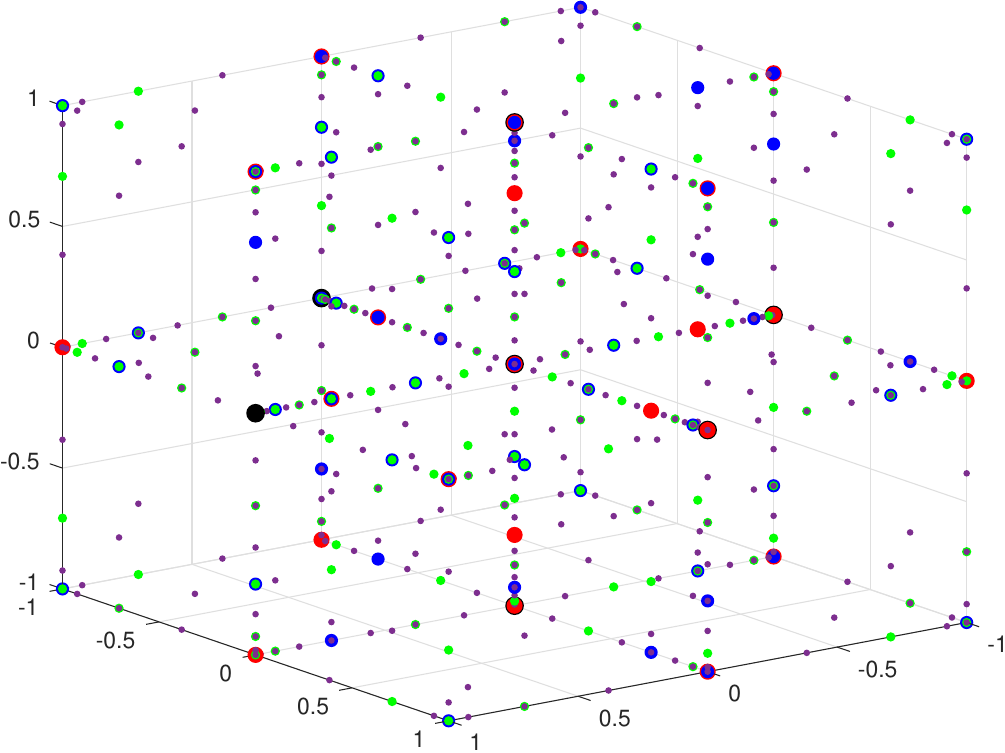}}
\end{minipage}}
\vskip -5pt
\caption
{\small  The sparse grid constructed from the Chebyshev-Gauss-Lobatto quadrature. Left: 2D sparse grid nodes; Right: 3D sparse grid nodes.}\label{sparse_grid_nodes}
\end{figure}

The Smolyak algorithm constructs the interpolant by combining all tensor-product interpolants that satisfy the condition $m - n + 1 \leq |\bm{l}| \leq m$.  Figure\,\ref{sparse_grid_nodes} shows a sparse grid based on the one-dimensional Chebyshev-Gauss-Lobatto quadrature in the interpolation space. This method not only enhances computational efficiency but also ensures the accuracy of the interpolation while avoiding the high costs associated with evaluating all possible points in high-dimensional spaces. Moreover, since it allows different resolutions in different dimensions, this approach provides the flexibility needed to address practical problems. Here, we develop a fully discrete Monte Carlo algorithm based on sparse grids for volumetric-constrained semilinear partial differential equations \eqref{mainprob}. To achieve this, we express the solution $u$ as the conditional expectation over each subinterval $[t_{i-1}, t_{i}]$ by employing the Feynman-Kac formula \eqref{FKi}. This representation not only highlights the probabilistic nature of the solution but also serves as the cornerstone of our computational framework. We then simulate the trajectories of the random process $X_{t}$ within each subinterval according to \eqref{Xi}. This involves generating paths that reflect the stochastic dynamics dictated by the underlying equation. By appropriately processing these simulated trajectories (specifically by averaging over multiple realizations), we derive an approximate solution for $u$ on the specified time grid $t_i$. At each time level, we use sparse grids to select spatial nodes for the calculation of the interpolation functions.

\subsection{The fully discrete scheme}
 Based on the random variable $\bz$ sampling methodology detailed in Section \ref{sec3.2}, the approximation of the trajectory $\{\tX_{t_i}\}$ in Section \ref{sec3.3} and the sparse grid interpolation outlined in Section \ref{sec3.1} , we are now well-positioned to develop the fully discrete scheme for the semi-linear PIDEs. 
To this end, we use the trapezoidal rule to approximate the integral of \eqref{FKi} and employ an implicit-explicit technique for the nonlinear term, namely, $u^i(t_{i}, \tX_{t_i})$ in \eqref{FKi} is replaced by $u^{i-1}_\ast(\tX_{t_i})$.
As a consequence,  the numerical solution of \eqref{FKi} (denoted by $u_\ast^i(\bx)$) at time grid $t_i$ reads
\begin{equation}\label{uapp1}
\begin{split}
u_\ast^{i}(\bx)= \;&
 \mathbb{E}_{\tX_{t_{i-1}}=\bx}\Big[ u_\ast^{i-1}(\tX_{t_{i}} ) \mathbb{I}_{\tau_{i}>t_{i}} +g(t_{i},\tX_{t_{i}} )\mathbb{I}_{\tau_{i}\leq t_{i}}\Big] \\&+ \mathbb{E}_{\tX_{t_{i-1}}=\bx}\Big[(\tau_{i}\land t_{i} -t_{i-1})f(u_\ast^{i-1}(\tX_{t_{i}}))\Big].
\end{split}
\end{equation}
\begin{remark}{\em 
Indeed, the above approximation scheme can be interpreted as a stochastic discretization of the following continuous equation
\begin{equation}
\label{bpi}
\begin{cases}
\partial_{t}u^i(t,\bx) - \mathcal{L}u^i(\bx,t)  =f(u^{i-1}(\bx,t)),\;&(\bx,t) \in\Omega\times\mathcal{T}_i,\\[4pt]
u(\bx,t_{i-1}) =u^{i-1}_\ast(\bx),\quad &\bx\in  \Omega,\\[4pt]
u(\bx,t) = g(\bx,t),\;&(\bx,t) \in\Omega^{c}\times[t_{i-1},t_{i}],
\end{cases}
\end{equation}
On the other hand, we recall the first order implicit-explicit scheme for time discretization of semilinear parabolic equations:
\begin{equation*}
\frac{u^{i}-u^{i-1}}{\Delta t}-\mathcal{L}u^{i}=f(u^{i-1}).
\end{equation*}
We observe that it shares a common characteristic with \eqref{bpi}, namely that the nonlinear term is treated explicitly, while the linear term is treated implicitly. This similarity underscores why the proposed method is termed as an implicit-explicit scheme. 
}
\end{remark}
The discrete scheme \eqref{uapp1} is fundamentally characterized by the explicit treatment of the nonlinear term, offering the advantage of avoiding iterative procedures typically required to address nonlinearities (see, e.g., \cite{yang2023}). Nevertheless, this approach requires prior knowledge of the numerical solution at the previous time step $t_{i-1}$ for any arbitrary location $\tX_{t_i}$, as $\tX_{t_i}$ is stochastically generated. However, in practice, the numerical solution at the previous time step $t_{i-1}$ is only available at the grid points, i.e., $\{u^{i-1}_\ast(\bx_\ell)\}$. 
Hence, at each time interval $\mathcal{T}_{i}$, it is essential to compute the numerical solution at the grid points $\{\bx_\ell\}_{\ell=1}^{N_s}$, after which the interpolation function based on the sparse grid $\Is[u^{i-1}_{\ast}]$ can be constructed. This interpolation function serves as the initial value function for the subsequent time step. By leveraging the properties of sparse grids, this approach significantly reduces computational complexity. Building on this setup, we can formulate a Monte Carlo procedure that utilizes random samples as
\begin{equation}
\label{Sii}
\begin{split}
S_{j}(t_{i},\bx)=& \mathcal{I}_S^d[u^{i-1}_{\ast}](\tX_{t_{i},j} )\mathbb{I}_{\tau_{i}>t_{i}}  +  g(t_{i},\tX_{t_{i},j})\mathbb{I}_{\tau_{i}\leq t_{i}} +(\tau_{i}\land t_{i} -t_{i-1})f(\mathcal{I}_S^d[u^{i-1}_{\ast}](\tX_{t_{i},j})).
\end{split}
 \end{equation}
where the index $j$ represents the $j$-th experiment. By applying the law of large numbers, the numerical solutions can be obtained by the following fully discretized scheme
  \begin{equation}\begin{split}
  \label{ui}
u_\ast^{i}(\bx) &=\lim_{M\rightarrow\infty}\frac{1}{M}\sum_{j=1}^{M}S_{j}(t_{i},\bx)
\\&=\lim_{M\rightarrow\infty}\frac{1}{M}\sum_{j=1}^{M}\Big[\Is[u^{i-1}_{\ast}](\tX_{t_{i},j} )\mathbb{I}_{\tau_{i}>t_{i}} \\&\;\;\quad
 +  g(t_{i},\tX_{t_{i},j})\mathbb{I}_{\tau_{i}\leq t_{i}} +(\tau_{i}\land t_{i} -t_{i-1})f(\Is[u^{i-1}_{\ast}](\tX_{t_{i},j}))\Big],
\end{split} \end{equation}
where $M$ represents the total number of experiments. This approach ensures the convergence to the exact solution by performing sufficient independent simulations of the stochastic process. As the number of simulated paths increases, the numerical results progressively approximate the expected value, thereby improving the accuracy and reliability of the numerical solutions. Notably, parallel computing can be utilized in practical implementations to enhance efficiency by simultaneously handling simulation paths and spatial grid points. Finally, we summarize the above derivations into the implicit-explicit Monte Carlo algorithm presented in {\bf Algorithm 1}.


 \begin{algorithm}
\caption{The implicit-explicit sparse Monte Carlo method for \eqref{mainprob}.}  \label{alg:Framwork}
\begin{algorithmic}[1]
\Require $\{\bx_\ell\}_{\ell=1}^{N_s}$: the set of sparse grid; $M$: the number of simulation paths; $\Delta t$: time step size.  \;\; 
\Ensure the numerical solution $\{u^N_\ast(\bx_\ell)\}_{\ell=1}^{N_s}$ at the terminal time.
\For{$i=1:N$}
\For{$\ell=1:N_{s}$ ({\bf in parallel})}
\For{$j=1:M$ ({\bf in parallel})}
  \State Set $\tX_{t_{i-1},j}:=\tX_{t_{i-1},j}^\ell=\bx_{\ell}$;
   \State Randomly generate a number $N_{\Delta t}$. 
\If{$\mathbb{P}(N_{\Delta t} = 0)$}
        \State Compute $\tX_{t_{i},j} = \bx_{\ell}  + \mu(t_{i-1}, \bx_\ell)\Delta t+ \sigma(t_{i-1}, \bx_\ell) W_{\Delta t} $;
\Else
        \State Simulate random variables $ \{\bz_k\}_{k=1}^{N_{\Delta t}} $ using Lemma \ref{lemker1}-\ref{lemker3}.
        \State Compute $\tX_{t_{i},j} = \bx_{\ell}  + \mu(t_{i-1}, \bx_\ell)\Delta t+ \sigma(t_{i-1},\bx_\ell) W_{\Delta t}+\sum_{k=1}^{N_{\Delta t}}c(t_{i-1},\bx_\ell,\bz_k)$;
\EndIf
\If{$\tX_{t_{i},j}\in \Omega$}

            \State Compute $S_{j}(t_{i},\bx_{\ell})=\Is[u^{i-1}_{\ast}](\tX_{t_{i},j} )+\Delta t f\big(\Is[u^{i-1}_{\ast}](\tX_{t_{i},j})\big)$;
\Else
            \State Compute $S_{j}(t_{i},\bx_{\ell})= g(t_{i},\tX_{t_{i,j}}) +\Delta t f\big(\Is[u^{i-1}_{\ast}](\tX_{t_{i},j})\big)$;
 \EndIf

\EndFor : in $j$
\State Compute $u_\ast^{i}(\bx_{\ell})\approx\frac{1}{M}\sum_{j=1}^{M}S_{j}(t_{i},\bx_{\ell})$;  
\EndFor : in $\ell$
  \State Using the numerical solution $\{ u_\ast^{i}(\bx_{\ell})\}_{\ell=1}^{N_{s}}$ to compute the interpolation function $\Is[u_\ast^{i}](\cdot)$; 
   \EndFor : in $i$  
\end{algorithmic}
\end{algorithm}

 \section{Error estimates.}
 In this section, we derive the rigorous error estimates for the proposed fully discrete explicit-implicit Monte Carlo method based on the sparse grid, where we focus on the bounded domain case with homogeneous global boundary conditions, i.e., $g(t,\bx) = 0$, as the method can be easily extended to inhomogeneous boundary conditions.  It should be noted that we do not consider the error arising from the stochastic discretization, i.e., the error between $\{\tX_{t_i}\}$ and $\{X_{t_i}\}$. In fact, incorporating this error into the analysis would make it excessively complex, far beyond the scope of this paper, and would distract from the main topic; please refer to a similar situation in \cite{yang2023}.
Furthermore, in the subsequent analysis, we assume that the nonlinear function $f(u)$ satisfies the following Lipschitz condition:
\begin{equation}\label{Lipscond}
|f(u_1)-f(u_2)| \leq L|u_1-u_2|, \quad L > 0,
\end{equation}
as well as the following two additional assumptions:

\noindent{\bf Assumption 1}: Let $\mu,\sigma,c \in C^{1}$, for any $t\in[0,T]$, assume that the drift coefficient $\mu$ satisfies a one-sided Lipschitz condition,  the that both diffusion coefficient $\sigma$ and jump coefficient $c$ satisfy global Lipschitz conditions: 
\begin{eqnarray}
\label{cond1}
 &&\langle \bx-\by, \mu(\cdot,\bx)- \mu(\cdot, \by)  \rangle  \leq C_{\mu}|\bx-\by|^{2}, \quad \text{for\; all}\;\; \bx,\by\in\R^{d},\\
 \label{cond2}
&&|\sigma(\cdot, \bx) - \sigma(\cdot, \by)|\leq C_{\sigma}|\bx-\by|^{2} , \quad \text{for\; all} \;\;\bx,\by\in\R^{d}, \\
\label{cond3}
&&|c(\cdot, \cdot,\bz_{1}) - c(\cdot, \cdot,\bz_{2})|\leq C_{c}|\bz_{1}-\bz_{2}|^{2} , \quad \text{for \; all} \;\;\bz_{1},\bz_{2}\in\R^{d}.
\end{eqnarray}
\noindent{\bf Assumption 2 (Uniform Ellipticity)}: 
There exist constants $0<\upsilon\leq\kappa<\infty$, such that for all $\bx\in \R^d$,
\begin{equation}
\label{cond4}
\upsilon^2\mathbb{I}_d\preceq \sigma(\bx,\cdot) \sigma(\bx,\cdot)^{\top} \preceq \kappa^2 \mathbb{I}_d \Longleftrightarrow \upsilon^2 |\bs{\xi}|^2\leq \bs{\xi}^\top (\sigma \sigma^\top) \bs{\xi} \leq \kappa^2 |\bs{\xi}|^2,
\end{equation}
for any $ \bs{\xi} \in \R^d \setminus \{\mathbf{0}\}$.

\medskip
 
Let $e^i(\bx)=u(\bx,t_i)-u^i_\ast(\bx)$, $i=0,\cdots,N$. Then, we have the following Lemma. 
\begin{lemma}\label{lmmdecom} {\em Let $u(t_i,\bx)$ and $u_\ast^i(\bx)$ be the solution of  \eqref{FKi} and \eqref{ui} at time grid $t_{i}$, respectively. 
Then it holds that
\begin{equation} 
\label{mainerr}
\|e^{i}\|_{\L2(\Omega)}= \|u(t_i,\cdot)-u_\ast^{i}\|_{\L2(\Omega)}\leq B^{i}_{1} + B^{i}_{2}+B^{i}_{3} +B^{i}_{4} + B^{i}_{5} +B^{i}_{6} ,
 \end{equation}
where 
\begin{equation}
\label{EI}\begin{split}
&B^{i}_{1}= \Big\|   \mathbb{E}_{\bx}\Big[ u^{i-1}(\tX_{t_{i}}) \mathbb{I}_{\tau_{i}>t_{i}}   -\Is[u^{i-1}_{\ast}](\tX_{t_{i}}) \mathbb{I}_{\tau_{i}>t_{i}}   \Big] \Big\|_{\L2(\Omega)}, \\
&B^{i}_{2}= \Big\|  \mathbb{E}_{\bx}\Big[\Is[u^{i-1}_{\ast}](\tX_{t_{i}}) \mathbb{I}_{\tau_{i}>t_{i}}  + g(\tau_{i},\tX_{\tau_{i}}) \mathbb{I}_{\tau_{i}\leq t_{i}}   \Big]   \\
&\quad- \frac{1}{M}\sum_{j=1}^{M}\Big[\Is[u^{i-1}_{\ast}](\tX_{t_{i},j} )\mathbb{I}_{\tau_{i}>t_{i}}  +  g(t_{i},\tX_{t_{i,j}})\mathbb{I}_{\tau_{i}\leq t_{i}}\Big]  \Big\|_{\L2(\Omega)}, \\
&B^{i}_{3}= \Big\| \mathbb{E}_{\bx}\Big[   \int_{t_{i-1}}^{t_{i} \land \tau_{i}} f( u(s,\tX_{s}))\,\d s  - \int_{t_{i-1}}^{t_{i}\land \tau_{i}}f( u^{i-1}(\tX_{s}))\,\d s   \Big] \Big\|_{\L2(\Omega)},\\
&B^{i}_{4}=\Big\| \mathbb{E}_{\bx}\Big[ \int_{t_{i-1}}^{t_{i}\land \tau_{i}}f( u^{i-1}(\tX_{s}))\,\d s  \Big]  -\mathbb{E}_{\bx}\Big[(\tau_{i}\land t_{i} - t_{i-1})f(u^{i-1}(\tX_{t_{i}}))\Big]   \Big\|_{\L2(\Omega)}, \\
&B^{i}_{5}= \Big\| \mathbb{E}_{\bx}\Big[(\tau_{i}\land t_{i} - t_{i-1})f(u^{i-1}(\tX_{t_{i}}))\Big] -\mathbb{E}_{\bx}\Big[(\tau_{i}\land t_{i} - t_{i-1})f(\Is[u^{i-1}_{\ast}](\tX_{t_{i}})) \Big] \Big\|_{\L2(\Omega)},\\
&B^{i}_{6}= \Big\|\mathbb{E}_{\bx}\Big[(\tau_{i}\land t_{i} - t_{i-1})f(\Is[u^{i-1}_{\ast}](\tX_{t_{i}}))\Big]-  \frac{1}{M}\sum_{j=1}^{M}(\tau_{i}\land t_{i} - t_{i-1})f(\Is[u^{i-1}_{\ast}](\tX_{t_{i},j}) ) \Big\|_{\L2(\Omega)},
\end{split}\end{equation}
with $\mathbb{E}_{\bx}$ shorthand for $\mathbb{E}_{\tX_{t_{i-1}}=\bx}$, and the initial $\tX_{t_{i-1},j}=\bx.$}
\end{lemma}
\vspace{-20pt}
\begin{proof}
By comparing the exact solution \eqref{FKi}, with $X_{t_i}$ replaced by $\tX_{t_i}$, and the numerical solution \eqref{ui} at the time grid $t_i$, then the error function is given by their difference:
\begin{eqnarray*}
&&\|e^{i}\|_{\L2(\Omega)}=\Big\|u^i(\cdot)  -  \frac{1}{M}\sum_{j=1}^{M}S_{j}(t_{i},\cdot)\Big\|_{\L2(\Omega)}\nonumber \\
&&= \Big\|   \mathbb{E}_{\bx}\Big[ u^{i-1}(t_{i-1},\tX_{t_{i}}) \mathbb{I}_{\tau_{i}>t_{i}}  + g(\tau_{i},\tX_{\tau_{i}}) \mathbb{I}_{\tau_{i}\leq t_{i}}    +  \int_{t_{i-1}}^{t_{i} \land \tau_{i}} f( u(s,\tX_{s}))\,\d s  \Big] \nonumber  \\
&&\quad- \frac{1}{M}\sum_{j=1}^{M}\Big[\Is[u^{i-1}_{\ast}](\tX_{t_{i},j} )\mathbb{I}_{\tau_{i}>t_{i}}  +  g(t_{i},\tX_{t_{i,j}})\mathbb{I}_{\tau_{i}\leq t_{i}} +(\tau_{i}\land t_{i} -t_{i-1})[f(\Is[u^{i-1}_{\ast}])] \Big]  \Big\|_{\L2(\Omega)}\\
&&\leq\Big\|\mathbb{E}_{\bx}\Big[ u^{i-1}(t_{i-1},\tX_{t_{i}}) \mathbb{I}_{\tau_{i}>t_{i}}  + g(\tau_{i},\tX_{\tau_{i}}) \mathbb{I}_{\tau_{i}\leq t_{i}}   \Big]   \nonumber\\
&&\quad- \frac{1}{M}\sum_{j=1}^{M}\Big[\Is[u^{i-1}_{\ast}](\tX_{t_{i},j} )\mathbb{I}_{\tau_{i}>t_{i}}  +  g(t_{i},\tX_{t_{i,j}})\mathbb{I}_{\tau_{i}\leq t_{i}}  \Big]  \Big\|_{\L2(\Omega)}  \nonumber  \\
 &&\quad+\Big\|  \mathbb{E}_{\bx}\Big[   \int_{t_{i-1}}^{t_{i} \land \tau_{i}} f(u(s,\tX_{s}))\,\d s  \Big]   -  \frac{1}{M}(\tau_{i}\land t_{i} -t_{i-1})\sum_{j=1}^{M}f(\Is[u^{i-1}_{\ast}])   \Big\|_{\L2(\Omega)}.\nonumber
\end{eqnarray*}
Next, we add the following three auxiliary terms in the above equation
\begin{equation*}\begin{split}
&\mathbb{E}_{\bx}\Big[ \int_{t_{i-1}}^{t_{i}\land \tau_{i}}f( u^{i-1}(\tX_{s}))\,\d s \Big] ,\;\;\;
\mathbb{E}_{\bx}\Big[(\tau_{i}\land t_{i} - t_{i-1})f(u^{i-1}(\tX_{t_{i}}))\Big],\\&\mathbb{E}_{\bx}\Big[(\tau_{i}\land t_{i} - t_{i-1})f(\Is[u^{i-1}_{\ast}](\tX_{t_{i}}))\Big],
\end{split}\end{equation*}
to obtain the desired result. This ends the proof.
\end{proof}

 Next, we estimate the bound of the six terms $\{B^{i}_j\}_{j=1}^6$ one by one as follows.
\begin{lemma}\label{lemB1} {\em
Let $u(t_{i-1},\bx)$ and $u_\ast^{i-1}(\bx)$ be the solution of \eqref{pi} and \eqref{ui} at time grid $t_{i-1}$, respectively.  Assume that $u_t\in L^\infty(I;L^2(\Omega))$, $\sigma$ satisfy the conditions of \eqref{cond2}and \eqref{cond4}, and $\kappa\ge (\Delta t)^{-1/2}$, then it holds that
\begin{equation}\label{boundB1}
B^{i}_{1}\leq  \| e^{i-1}\|_{\L2(\Omega)}. 
\end{equation}}
\end{lemma}
\vspace{-20pt} 
\begin{proof}
Since the expectation term in $B^{i}_{1}$ is conditioned on $\tX_{t_{i-1}}=\bx$ and the process $\tX_{t}$ has independent and stationary increment, it follows that $\tX_{t_{i-1}}| \tX_{t_{i-2}}=\bx$ and $\tX_{t_{i}}| \tX_{t_{i-1}}=\bx$ have the same distribution for fixed time step size $\Delta t$. Hence, we find from \eqref{EI} that 
\begin{eqnarray}\label{B1proof1}
&& B^{i}_{1}= \Big\|   \mathbb{E}_{\tX_{t_{i-1}}=\bx}\Big[ u^{i-1}(t_{i-1},\tX_{t_{i}}) \mathbb{I}_{\tau_{i}>t_{i}}   -\Is[u^{i-1}_{\ast}](\tX_{t_{i}}) \mathbb{I}_{\tau_{i}>t_{i}}   \Big] \Big\|_{\L2(\Omega)}\nonumber\\
 &&\hspace{15pt}= \Big\|   \mathbb{E}_{\tX_{t_{i-2}}=\bx}\Big[ u^{i-1}(t_{i-1},\tX_{t_{i-1}}) \mathbb{I}_{\tau_{i-1}>t_{i-1}}   -\Is[u^{i-1}_{\ast}](\tX_{t_{i-1}}) \mathbb{I}_{\tau_{i-1}>t_{i-1}}   \Big] \Big\|_{\L2(\Omega)}\nonumber\\
 &&\hspace{15pt}=  \Big\|   \mathbb{E}_{\tX_{t_{i-2}}=\bx}\Big[ u^{i-1}(t_{i-1},\tX_{t_{i-1}})  -\Is[u^{i-1}_{\ast}](\tX_{t_{i-1}})\Big |  \tau_{i}>t_{i}\Big]  \Big\|_{\L2(\Omega)}\nonumber\\
 &&\hspace{15pt}\leq \Big\|   \mathbb{E}_{\tX_{t_{i-2}}=\bx}\Big[ u^{i-1}(t_{i-1},\tX_{t_{i-1}})  -\Is[u^{i-1}_{\ast}](\tX_{t_{i-1}}) \Big]  \Big\|_{\L2(\Omega)}\\
 &&\hspace{15pt}=\Big\|  \mathbb{P}(N_{\Delta t}=0) \mathbb{E}_{\tX_{t_{i-2}}=\bx}\Big[ u^{i-1}(t_{i-1},\tX_{t_{i-1}})  -\Is[u^{i-1}_{\ast}](\tX_{t_{i-1}}) \Big]  \nonumber\\
 &&\hspace{15pt}\quad+  \sum_{k=1}^{\infty}\mathbb{P}(N_{\Delta t}=k) \mathbb{E}_{\tX_{t_{i-2}}=\bx}\Big[ u^{i-1}(t_{i-1},\tX_{t_{i-1}})  -\Is[u^{i-1}_{\ast}](\tX_{t_{i-1}}) \Big]  \Big\|_{\L2(\Omega)} \nonumber\\
 &&\hspace{15pt}\leq B^{i}_{1,1}+B^{i}_{1,2},\nonumber
\end{eqnarray}
where 
\begin{equation}\label{B11B12} 
\begin{split}
&B^{i}_{1,1}=\Big\|  \mathbb{P}(N_{\Delta t}=0) \mathbb{E}_{\tX_{t_{i-2}}=\bx}\Big[ u^{i-1}(t_{i-1},\tX_{t_{i-1}})  -\Is[u^{i-1}_{\ast}](\tX_{t_{i-1}}) \big| N_{\Delta t}=0\Big]  \Big\|_{\L2(\Omega)} , \\
    &B^{i}_{1,2}=\Big\|  \sum_{k=1}^{\infty}\mathbb{P}(N_{\Delta t}=k)  \mathbb{E}_{\tX_{t_{i-2}}=\bx}\Big[ u^{i-1}(t_{i-1},\tX_{t_{i-1}})  -\Is[u^{i-1}_{\ast}](\tX_{t_{i-1}})  \big| N_{\Delta t}=k\Big]  \Big\|_{\L2(\Omega)}.
    \end{split}
\end{equation}

We first prove the bound for $B^{i}_{1,1}$ with the case $N_{\Delta t} = 0$, in which the process is driven by Brownian motion without jumps. Thus, the governing equation does not involve any nonlocal operators. Under the assumptions of uniform ellipticity and bounded drift, the transition probability density $p_{\Delta t}(\bx, \by)$ satisfies the following estimate (cf. \cite[p.\,895]{Aronson1967})
\begin{equation*}
   p_{\Delta t}(\bx,\by) \leq  C_e(\Delta t)^{-d/2} e^{-\frac{|\by-\bx|^2}{C_d\Delta t}},
\end{equation*}
where the positive constants $C_e,C_d$ depend on the ellipticity constant and the upper bound of the drift term. From the Gaussian-type upper bound and the substitution $\xi = \bx - \by$, we obtain that

\begin{equation}\label{intdensity}
\begin{split}
\int_{\R^d} \Big( p_{\Delta t}(\bx,\by) \Big)^{2}\,\d\by 
&\leq   C_e^2(\Delta t)^{-d} \int_{\R^{d}} e^{-\frac{2|\bx-\by|^2}{C_d \Delta t}}\,\d\by  \\
&\leq C_e^2(\Delta t)^{-d} \int_{\R^{d}} e^{-\frac{2\bs{\xi}^{2}}{C_d \Delta t}}\,\d\bs{\xi} 
\\&\leq C_e^2(\Delta t)^{-d} \Big(\frac{C_d \Delta t}{2}\Big)^{d/2} \int_{\R^{d}} e^{-\bs{\eta}^{2}}\,\d\bs{\eta}
\\&=C_e^2\Big(\frac{\pi C_d}{2\Delta t}\Big)^{d/2}.
\end{split}
\end{equation}
According to the assumptions \eqref{cond4}, if the ellipticity constant $\kappa$ satisfies $\kappa \geq (\Delta t)^{-1/2}$, then by choosing the constants as $C_d = 4/\kappa,$ $C_e^2 = \frac{1}{(2\pi \kappa)^{d/2}},$  we obtain the estimate $\int_{\mathbb{R}^{n}} \big(p_{\Delta t}(\bx,\by)\big)^{2}\,d\by \leq 1.$
Thus, we derive from Cauchy-Schwarz inequality, \eqref{B11B12}, and \eqref{intdensity} that
\begin{equation}
\label{EuQ}
\begin{split}
& \mathbb{E}_{\tX_{t_{i-2}}=\bx}\Big[ u^{i-1}(t_{i-1},\tX_{t_{i-1}})  -\Is[u^{i-1}_{\ast}](\tX_{t_{i-1}})  \big| N_{\Delta t}=0\Big]  \\
&= \int_{\R^{d}} \big( u^{i-1}(t_{i-1},\by)  -\Is[u^{i-1}_{\ast}](\by) \big)p_{\Delta t}(\bx,\by)\,\d\by \\
 & \leq  \Big(  \int_{\R^{d}} \big( u^{i-1}(t_{i-1},\by)  -\Is[u^{i-1}_{\ast}](\by) \big)^{2}\,\d\by   \Big)^{1/2} \cdot\Big(  \int_{\R^{d}} \big( p_{\Delta t}(\bx,\by) \big)^{2}\,\d\by  \Big)^{1/2}
 \\&\leq \Big(  \int_{\R^{d}} \big( u^{i-1}(t_{i-1},\by)  -\Is[u^{i-1}_{\ast}](\by) \big)^{2}\,\d\by   \Big)^{1/2}.
\end{split}
\end{equation}
Thanks to \eqref{B11B12}, \eqref{EuQ}, and  the fact that $\mathbb{P}(N_{\Delta t}=0) = 1-e^{-\lambda \Delta t}<1$, we find that
\begin{equation}
\label{B11}
\begin{split}
B^{i}_{1,1} 
&\leq \Big(\int_{\Omega}  \mathbb{E}_{\bx}\Big[ u^{i-1}(t_{i-1},\tX_{t_{i-1}})  -\Is[u^{i-1}_{\ast}](\tX_{t_{i-1}}) \Big]^2  \,\d\bx\Big)^{1/2}\\
&\leq\Big(\int_{\Omega}  \Big(  \int_{\R^{d}} \big( u^{i-1}(t_{i-1},\by)  -\Is[u^{i-1}_{\ast}](\by) \big)^{2}\,\d\by \Big)  
\Big)^{1/2} 
\leq  \| e^{i-1}\|_{L^2(\Omega)}.
\end{split}
\end{equation}
Moreover, it is easy to verify that the probability of the Poisson process $N_{\Delta t}$ having $k$ jumps within the interval $[t_{i-1}, t_{i}]$, as given in \eqref{PNk}, behaves as 
$$\sum_{k=1}^{\infty} \mathbb{P}(N_{\Delta t} = k) = \mathcal{O}(\Delta t),$$
which, together with \eqref{B11B12}, implies that
\begin{equation}
\label{B12}
\begin{split}
B^{i}_{1,2}&=\Big\| \sum_{k=1}^{\infty}\mathbb{P}(N_{\Delta t}=k)  \mathbb{E}_{\bx}\Big[ u^{i-1}(t_{i-1},\tX_{t_{i-1}})  -\Is[u^{i-1}_{\ast}](\tX_{t_{i-1}})  \big| N_{\Delta t}=k\Big]  \Big\|_{L^2(\Omega)} \\
&\leq \sum_{k=1}^{\infty}\mathbb{P}(N_{\Delta t}=k)  \Big\| \mathbb{E}_{\bx}\Big[ u^{i-1}(t_{i-1},\tX_{t_{i-1}})  -\Is[u^{i-1}_{\ast}](\tX_{t_{i-1}})  \big| N_{\Delta t}=k\Big]  \Big\|_{L^2(\Omega)}
\\&\leq \Delta t  \| e^{i-1}\|_{L^2(\Omega)}.
\end{split}
\end{equation}
A combination of the above results leads to the desired result.
\end{proof}
 
\begin{lemma}
\label{lemB3} {\em
Let $u(\bx,t_{i-1})$ be the solution of \eqref{mainprob} at time grid $t_{i-1}$.  
We assume that $u_{t}\in L^\infty(\mathcal{T}_i;L^2(\Omega))$ and $f(\cdot)$ fulfills the Lipschitz condition \eqref{Lipscond}. Then, it holds that

 \begin{equation}
 \label{bounB3}
 B^{i}_{3}\leq C (\Delta t)^2\|u_t\|_{L^\infty(\mathcal{T}_i;L^2(\Omega))},
 \end{equation}
where $C$ is a positive constant independent of $\Delta t$ and $u$.}
\end{lemma}
\vspace{-20pt}
\begin{proof}
Since the expectation is taken in $ L^2(\Omega) $, we can leverage linearity and additivity to rewrite $ B^i_3 $ as:
\[
B^{i}_3 = \Big\| \mathbb{E}_{\tX_{t_{i-1}}=\bx}\Big[ \int_{t_{i-1}}^{t_{i} \land \tau_{i}} \Big( f(  u(s, \tX_{s})) - f( u^{i-1}(t_{i}, \tX_{s})) \Big) \, \mathrm{d}s \Big] \Big\|_{L^2(\Omega)}.
\]
Using the mean-value Theorem and Cauchy-Schwarz inequality, we find that
\begin{eqnarray}
\label{Schwarzineq}
&&B^{i}_{3}
\leq\Big\| \mathbb{E}_{\tX_{t_{i-1}}=\bx}\Big[ \int_{t_{i-1}}^{t_{i} \land \tau_{i}} \Big( f(  u(s, \tX_{s})) - f( u^{i-1}(t_{i}, \tX_{s})) \Big) \, \mathrm{d}s \Big] \Big\|_{L^2(\Omega)}.\nonumber\\
&&\quad\leq L\Big\| \mathbb{E}_{\tX_{t_{i-1}}=\bx}\Big[  \int_{t_{i-1}}^{t_{i}\land \tau_{i}} \big|u(s, \tX_{s}) - u^{i-1}(t_{i}, \tX_{s})|\,\d s\Big]\Big\| _{L^2(\Omega)} \nonumber\\
&&\quad\leq L\|u_t\|_{L^\infty(\mathcal{T}_i;L^2(\Omega))}\Big\| \mathbb{E}_{\tX_{t_{i-1}}=\bx}\Big[ \int_{t_{i-1}}^{t_{i}\land \tau_{i}}(s-t_{i-1})\, \d s\Big]\Big\| _{L^2(\Omega)} \\
&&\quad\leq L\|u_t\|_{L^\infty(\mathcal{T}_i;L^2(\Omega))}\Big\| \mathbb{E}\Big[ (t_{i}\land \tau_{i}-t_{i-1})^2\Big]\Big\| _{L^2(\Omega)}\nonumber
\\&&\quad\leq C_{L}(\Delta t)^2\|u_t\|_{L^\infty(\mathcal{T}_i;L^2(\Omega))},\nonumber
\end{eqnarray}
which yields the error bound.
\end{proof}

 \begin{lemma}\label{lemB4} {\em 
Assume that $f(\cdot)$ satisfies the Lipschitz condition \eqref{Lipscond}, and $\mu,\sigma,c$ satisfy the conditions \eqref{cond1}-\eqref{cond3}, $\nabla^2 u\in L^\infty(I;\Omega)$, then we have
\begin{equation}
\label{boundB4}
B^{i}_{4}\leq C(\Delta t )^{2},
\end{equation}
where constant $C$ depends on upper bounds of $\mu$, $\sigma$, $c$, $f$ and their derivatives of $u$.}
\end{lemma}
\vspace{-20pt}
\begin{proof}
We derive from \eqref{EI}  and the Lipschitz condition \eqref{Lipscond} that
\begin{equation}\label{B4err}
\begin{split}
B^{i}_{4}&=\Big\| \mathbb{E}_{\bx}\Big[ \int_{t_{i-1}}^{t_{i}\land \tau_{i}}f( u^{i-1}(\tX_{s}))\,\d s  \Big]  -\mathbb{E}[\tau_{i}\land t_{i} - t_{i-1}]f(u^{i-1}(\tX_{t_{i}}))   \Big\|_{L^2(\Omega)}\\
&  \leq \Big\| \mathbb{E}_{\bx}\Big[ \int_{t_{i-1}}^{t_i \land \tau_i} \Big|f(u^{i-1}(\tX_s)) -f(u^{i-1}(\tX_{t_{i}}))  \Big|\,\d s  \Big] \Big\|_{L^2(\Omega)}
\\&\leq L \Big\| \mathbb{E}_{\bx}\Big[ \int_{t_{i-1}}^{t_i \land \tau_i} \big|u^{i-1}(\tX_s) -u^{i-1}(\tX_{t_{i}})  \big|\,\d s  \Big] \Big\|_{L^2(\Omega)}.
\end{split}
\end{equation}
For any $u^{i-1} \in C^2$, the variations of $u^{i-1}(\tX_s)$ can be expanded as 
\begin{equation}
\label{deltaui1}
u^{i-1}(\tX_{t_{i}}) = u^{i-1}(\tX_{s}) + \nabla u^{i-1}(\tX_{s}) \cdot (\tX_{t_{i}} -\tX_{s}) + \frac{1}{2} (\tX_{t_{i}} - \tX_{s})^\top \nabla^2 u^{i-1}(\xi) (\tX_{t_{i}} -\tX_{s}),
\end{equation}
where $ \xi \in [\tX_s, \tX_{t_i} ]$. On the other hand, the dynamics of stochastic process $\tX_t$ satisfies 
\begin{equation}
\label{XtiXs}
\tX_{t_i} -\tX_s =  \mu(s , \tX_{s} ) (t_{i}-s) + \sigma(s, \tX_{s} ) (W_{s} -W_{t_{i}}) + \sum_{k=1}^{N_{\Delta t}}c(t,\tX_{s},\bz_{k}).
\end{equation}
Since $\mu$ and $\sigma$ are assumed to be bounded, the increment of $\mu(s, \tX_{s})$ over a small time interval $[s,t_{i}]$ is of order $\mathcal{O}(t_{i} - s)$, and $\sigma(s, \tX_{s})(W_{t_{i}} - W_{s})$, where $(W_{t_{i}} - W_{s}) \sim \mathcal{N}(0, t_{i} - s)$, is of order $\mathcal{O}(\sqrt{t_{i} - s})$. Similarly, it is readily verified that the term $\sum_{k=1}^{N_{\Delta t}} c(t, X_{s}, \bz_{k})$ is of order $\mathcal{O}(t_{i} - s)$.

Substituting \eqref{XtiXs} into \eqref{deltaui1} and combining it with \eqref{B4err}, we obtain that
\begin{equation*}
\begin{split}
B_{4}^{i}  
& \leq L\Big\| \mathbb{E}_{\bx}\Big[ \int_{t_{i-1}}^{t_i \land \tau_i} \Big| \nabla u^{i-1}(\tX_{s}) \cdot (\tX_{t_{i}} - \tX_{s}) + \frac{1}{2} (\tX_{t_{i}} - \tX_{s})^\top \nabla^2 u^{i-1}(\xi) (\tX_{t_{i}} - \tX_{s}) \Big|\,\d s  \Big] \Big\|_{L^2(\Omega)} \\
&\leq L\Big\| \mathbb{E}_{\bx}\Big[ \int_{t_{i-1}}^{t_i \land \tau_i} | \nabla u^{i-1}(\tX_{s}) \cdot (\tX_{t_{i}} - \tX_{s})  |\,\d s  \Big] \Big\|_{L^2(\Omega)} \\
&\quad +L\Big\| \mathbb{E}_{\bx}\Big[ \int_{t_{i-1}}^{t_i \land \tau_i} \Big|\frac{1}{2} (\tX_{t_{i}} - \tX_{s})^\top \nabla^2 u^{i-1}(\xi) (\tX_{t_{i}} - \tX_{s}) \Big|\,\d s  \Big] \Big\|_{L^2(\Omega)} \\
&\leq L \Big\| \mathbb{E}_{\bx}\Big[ \int_{t_{i-1}}^{t_i \land \tau_i} |\nabla u^{i-1}(\tX_{s}) |\cdot\Big|  \mu(s ,  \tX_{s} ) (t_{i}-s) 
\\&\quad\quad+ \sigma(s,  \tX_{s} ) (W_{s} -W_{t_{i}}) + \sum_{k=1}^{N_{\Delta t}}c(t,\tX_{s},\bz_{k})\Big|\,\d s  \Big] \Big\|_{L^2(\Omega)} \\
&\quad +L\Big\| \mathbb{E}_{\bx}\Big[ \int_{t_{i-1}}^{t_i \land \tau_i} \Big|\frac{1}{2} (\tX_{t_{i}} -  X_{s})^\top \nabla^2 u^{i-1}(\xi) (\tX_{t_{i}} -  X_{s}) \Big|\,\d s  \Big] \Big\|_{L^2(\Omega)}.
\end{split}
\end{equation*}
Due to the martingale property of the Brownian motion, we have 
$$\mathbb{E}_{\bx}\Big[ \int_{t_{i-1}}^{t_i \land \tau_i}  \sigma(s,   X_{s} ) (W_{s} -W_{t_{i}}) \d s \Big] =0.$$
Moreover, the probability of the Poisson process $N_t$ experiencing $k$ jumps is characterized by \eqref{PNk}, and its upper bound is given by
$$\sum_{k=1}^{\infty} \mathbb{P}(N_{\Delta t} = k) = \mathcal{O}(\Delta t),$$ which implies 
$$L\Big\| \mathbb{E}_{\bx}\Big[ \int_{t_{i-1}}^{t_i \land \tau_i}|\nabla u^{i-1}( X_{s}) |\cdot \Big|  \mu(s ,   X_{s} ) (t_{i}-s) + \sum_{k=1}^{N_{\Delta t}}c(t, X_{s},\bz_{k})\Big|\,\d s  \Big] \Big\|_{L^2(\Omega)} \leq C(\Delta t)^{2}.$$
Similarly, we have  
\begin{equation}
    \mathbb{E}_{\bx}\Big[ \int_{t_{i-1}}^{t_i \land \tau_i} \Big|\frac{1}{2} ( X_{t_{i}} -  X_{s})^\top \nabla^2 u^{i-1}(\xi) ( X_{t_{i}} -  X_{s}) \Big|\,\d s  \Big] \leq C(\Delta t)^2,
\end{equation}
which establishes the error bound.
\end{proof}

\begin{lemma}
\label{lemB5} {\em
 Let $u(\bx,t_{i-1})$ and $u_\ast^{i-1}(\bx)$ be the solution of \eqref{mainprob} and \eqref{ui} at time grid $t_{i-1}$, respectively.
Assume that $f(\cdot)$ satisfies the Lipschitz condition \eqref{Lipscond}, $\sigma$ satisfy the conditions of \eqref{cond2}and \eqref{cond4}, and $\kappa\ge (\Delta t)^{-1/2}$, then  
\begin{equation}
\label{boundB5}
\begin{split}
B^{i}_{5} \leq   C\Delta t\|e^{i-1}\|_{L^2(\Omega)},
\end{split}
\end{equation}
where the positive constant $C$ is indedepent of $u$ and $\Delta t$.}
\end{lemma}\vspace{-20pt}
\begin{proof}
We derive from the Lipschitz condition, Cauchy-Schwarz inequality, \eqref{intdensity}, and \eqref{EI}  that
\begin{equation*}
\begin{split}
B^{i}_{5}
&=\Big\| \mathbb{E}_{\tX_{t_{i-1}}=\bx}\Big[(\tau_{i}\land t_{i} - t_{i-1})\big( f(u^{i-1}(\tX_{t_{i}}))  -f(\Is[u^{i-1}_{\ast}](\tX_{t_{i}})) \big) \Big] \Big\|_{L^2(\Omega)}
\\&=\Big(\int_{\Omega}\mathbb{E}_{\tX_{t_{i-1}}=\bx}\Big[(\tau_{i}\land t_{i} - t_{i-1})\big( f(u^{i-1}(\tX_{t_{i}}))  -f(\Is[u^{i-1}_{\ast}](\tX_{t_{i}})) \big) \Big]^2\,\d\bx\Big)^{\frac12} 
\\&\leq C \Delta t \Big(\int_{\Omega}  \Big(  \int_{\R^{d}} \big( u^{i-1}(t_{i-1},\by)  -\Is[u^{i-1}_{\ast}](\by) \big)^{2}\,\d\by \Big)  \cdot\Big(  \int_{\R^{d}} \Big( p_{\Delta t}(\bx,\by) \Big)^{2}\,\d\by  \Big)  \d\bx \Big)^{1/2}
\\&\leq C \Delta t\Big(\int_{\Omega} \int_{\R^{d}} \big( u^{i-1}(t_{i-1},\by)  -\Is[u^{i-1}_{\ast}](\by) \big)^{2} \,\d\by\d\bx \Big)^{1/2} 
\\&\leq C \Delta t\|e^{i-1}\|_{L^2(\Omega)}. 
\end{split}
\end{equation*}
 This ends the proof.
\end{proof}

With the above analysis, the error bound $\|e^N\|_{L^{2}(\Omega)}$ can be easily derived as below.
\begin{theorem}
\label{EN}
Let $u(\bx,t_i)$ and $u^i_\ast(\bx)$ be the solutions of \eqref{mainprob} and \eqref{ui}, respectively. Assume that $\nabla^2 u\in L^\infty(I;\Omega),$ $u_t\in L^\infty(I;L^2(\Omega))$, $f(\cdot)$ fulfills the Lipschitz condition \eqref{Lipscond} and $\mu,\sigma,c$ satisfy the conditions \eqref{cond1}-\eqref{cond4},  and $\kappa\ge (\Delta t)^{-1/2}$, then there holds 
\begin{equation}\begin{split}
\label{eN}
\|u(\cdot,T)-u^N_\ast\|_{L^2(\Omega)}\leq  CM^{-1/2} +C\Delta t\|u_t\|_{L^\infty(0,T;L^2(\Omega))}+C\Delta t,
\end{split}\end{equation}
where $C$ is a positive constant independent of $\Delta t$ and $M$.
\end{theorem}\vspace{-20pt}
\begin{proof}
For $g(\bx,t)=0$, we only need to consider the case with $\tau_{i}>t_{i}$, it is obvious that 
\begin{equation*}
\begin{split}
B^{i}_{2}  \leq\Big\|  \mathbb{E}_{\bx}\Big[\Is[u^{i-1}_{\ast}](\tX_{t_{i}}) \mathbb{I}_{\tau_{i}>t_{i}}    \Big] - \frac{1}{M}\sum_{j=1}^{M}\Big[\Is[u^{i-1}_{\ast}](\tX_{t_{i},j} )\mathbb{I}_{\tau_{i}>t_{i}}   \Big]  \Big\|_{L^2(\Omega)}\leq  CM^{-1/2}. 
\end{split}
\end{equation*}
Similarly, we have
\begin{equation*}
\begin{split}
B^{i}_{6}= \Big\|\mathbb{E}[\tau_{i}\land t_{i} - t_{i-1}]f(\Is[u^{i-1}_{\ast}](\tX_{t_{i}}))-  \frac{1}{M}\sum_{j=1}^{M}(\tau_{i}\land t_{i} - t_{i-1})f(\Is[u^{i-1}_{\ast}])(\tX_{t_{i},j})  \Big\|_{L^2(\Omega)} \leq CM^{-1/2}.
\end{split}
\end{equation*}
By \eqref{mainerr}, and Lemma \ref{lemB1}-Lemma\ref{lemB5}, we deduce that
\begin{equation}\label{eijk}
\begin{split}
\|e^{i}\|_{L^2(\Omega)} \leq  \| e^{i-1}\|_{L^2(\Omega)} +C M^{-1/2}+C (\Delta t)^2\|u_t\|_{L^\infty(\mathcal{T}_i;L^2(\Omega))}+C(\Delta t )^{2}.
\end{split}
\end{equation}
Summing \eqref{eijk} from $i=1$ to $i=N$ and using the fact that $\|e^{0}\|_{L^2(\Omega)}=0$, we can easily obtain the desired results \eqref{eN}. 
\end{proof}
\section{Numerical results}
In this section, we conduct a series of numerical results to rigorously validate the accuracy, efficiency, and stability of our fully discrete Monte Carlo algorithm for solving semilinear parabolic partial differential equations with volume constraints on bounded domains. To obtain exact solutions for benchmarking, we first consider a model with a constant kernel function, which allows us to verify the accuracy and convergence of the proposed method. Next, we present a 100-dimensional test case without jump terms to demonstrate that the algorithm maintains high numerical accuracy and stability in high-dimensional settings, thereby confirming its effectiveness for real-world problems involving large dimensions. Finally, we examine a 100-dimensional problem with a hypersingular integral kernel to demonstrate the method's flexibility and broad applicability.

All numerical results are produced by Algorithm~\ref{alg:Framwork}, which computes the numerical solution $u^N_{\ast}$ at a prescribed set of grid points $T$, along with the corresponding $L^2$-errors that measure the deviation from the exact solution:
\begin{equation}
\label{error}
\|e^N\|_{L^2(\Omega)}=\|u(\cdot,T)-u_\ast^N\|_{L^2(\Omega)}=\Big(\int_\Omega(u(\cdot,T)-u_\ast^N)^2\,\d \bx\Big)^{1/2}.
\end{equation}
 \begin{exa} \label{example1}{\bf (Accuracy test with constant kernel function)} We begin by considering problem \eqref{mainprob} on a high-dimensional cube $\Omega =(-1,1)^d$. The jump amplitude function is chosen as $c(t, \bx, \bz) = \bz$, with a constant kernel function $\varphi(\bz) = \mathbf{1}_{|\bz| \leq \delta}$, and the nonlinear forcing term is defined as $f(u) = \frac{1 - u}{1 + u}$. The drift and diffusion coefficients are specified as follows:
\begin{equation*}
 \begin{split}
&\mu(t,\bx) =\big(\mu(t,x_{1}),\cdots, \mu(t,x_{d})\big)^{\top},\;\;   \mu(t,x_{i}) =  \frac{1}{10}\cos(2\pi t)(2x_{i}^{7} - x_{i}^{8}),\\
&\sigma(t,\bx) = \frac{1}{20}e^{-t}\begin{pmatrix}
  \sin^{2}(\pi x_{1}) &\cdots&0 \\   
 \vdots&  \ddots& \vdots\\[6pt]
 0&\cdots& \sin^{2}(\pi x_{d}) 
\end{pmatrix}.
 \end{split}
 \end{equation*}
Under this setting, the problem admits the exact solution
 \begin{equation}\label{exactsolu1}
 u(t,\bx) =  \sin(10t)\Big(\sum_{j=1}^d\Big(\frac{1}{j}x_j^5-\frac{1}{j+1}x^3_j\Big)\Big),\;\;\;(t,\bx)\in [0,\infty)\times\Omega.
 \end{equation}
In fact, by substituting the exact solution $u(t,\bx)$ into the equation \eqref{mainprob}, we obtain a right-hand side source function $r(t,\bx)$, which is incorporated into the nonlinear term $f(u)$ for consistency in the formulation.
The interaction domain $\Omega_{E}$ is defined by the extension from $\Omega$ by a radius of the horizon $\delta$, where $\mathcal{D}=\{ \bz \in \R^{d} \big| | \bz| \leq \delta \} $. We set the terminal time $T=1$, $\delta =0.4$ and solve \eqref{mainprob} numerically. 
\end{exa}

We first test the convergence rate of the proposed algorithm with respect to both time discretization and Monte Carlo sampling. The spatial discretization error can be considered negligible compared to the time discretization and sampling errors, owing to the spectral accuracy of the approximation error in the spatial direction achieved through interpolation reconstruction using Chebyshev-Gauss-Lobatto quadrature-based sparse grids. For convenience, in the simulation, we use a sparse grid with level $\bm{l}=3$ to select spatial nodes, and apply the sparse grid interpolation rule \eqref{sgqr} to construct interpolations of the numerical solutions obtained at each time step.

\begin{figure}[!h]\hspace{0pt}
\subfigure{
\begin{minipage}[t]{0.4\textwidth}
\centering
\rotatebox[origin=cc]{-0}{\includegraphics[width=1.2\textwidth,height=1\textwidth]{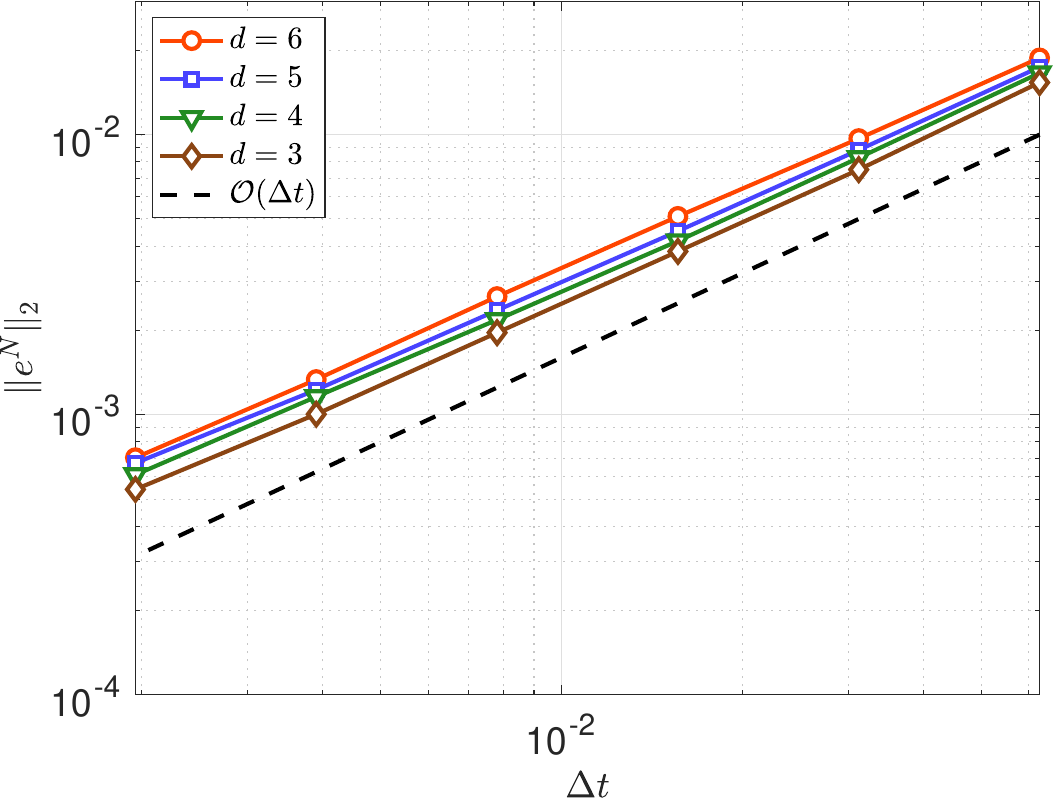}}
\end{minipage}}\hspace{30pt}
\subfigure{
\begin{minipage}[t]{0.4\textwidth}
\centering
\rotatebox[origin=cc]{-0}{\includegraphics[width=1.2\textwidth,height=1\textwidth]{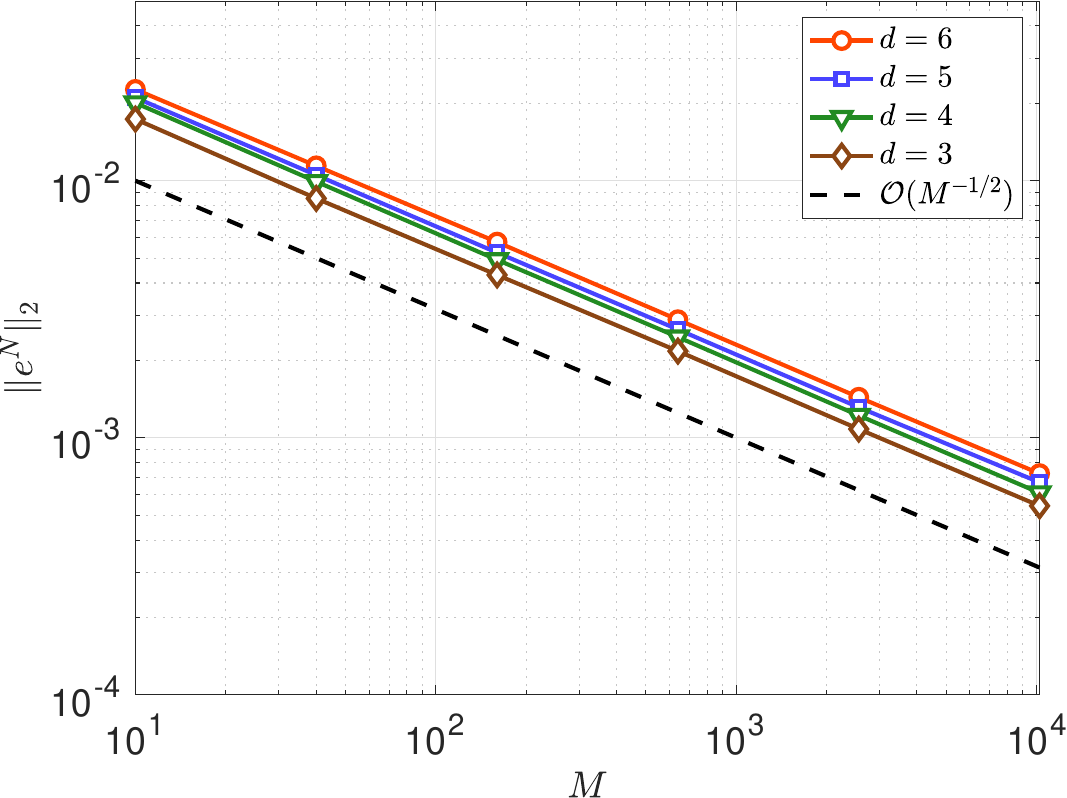}}
\end{minipage}}
\vskip -5pt
\caption
{\small  Left:  The numerical errors of \eqref{exactsolu1} against various time step $\Delta t$ with $T=1$ and $M=10000$; Right: The numerical errors of \eqref{exactsolu1} against the number of the path $M$ with $T=1$ and $\Delta t=2^{-10}$.}\label{err_exact}
\end{figure}

In Figure \ref{err_exact} (left), we show the discrete $L^2$-error of the solution to \eqref{exactsolu1} for various time step sizes $\Delta t$ and spatial dimensions $d = 3, 4, 5, 6$. We observe that as the time step $\Delta t$ decreases, the numerical error also decreases, with the convergence rate matching the estimate in \eqref{eijk}, approximately of order $\mathcal{O}(\Delta t)$. In Figure \ref{err_exact} (right), we present the discrete $L^2$-error for different numbers of simulation paths $M$ and dimensions $d = 3, 4, 5, 6$. As the number of paths $M$ increases, the numerical error decreases, showing half-order convergence, which is consistent with the analytical results discussed earlier. Both plots also show that as the spatial dimension $d$ increases, the overall computational error tends to increase.  

To demonstrate that the proposed method does not require satisfaction of the conventional time-space grid ratio condition, i.e., $(\Delta x)^2 \geq \Delta t$, we consider a two-dimensional example. A uniform spatial grid with fixed mesh size $\Delta x = 0.05$ is used. For several time step sizes $\Delta t = 2^{-3}, \ldots, 2^{-7}$, this condition is clearly violated. Nevertheless, the numerical results reported in Table~\ref{tab1} confirm that first-order accuracy is still achieved.

 \begin{table}[!htp]
\caption{The $L^2$-errors for various time step sizes $\Delta t$, with fixed spatial mesh size $\Delta x = 0.05$, dimension $d = 2$, and number of simulation paths $M = 10000$.}
\setlength{\tabcolsep}{5mm}{
\begin{tabular}{ c c c c c c }
\hline 
 & $\Delta t = 2^{-3}$& $\Delta t = 2^{-4}$& $\Delta t = 2^{-5}$& $\Delta t = 2^{-6}$& $\Delta t = 2^{-7}$\\[3pt]
\hline
$L^{2}$-error& 4.5216e-2 &2.2814e-2&1.1746e-2 & 6.0830e-3 & 3.1090e-3\\[3pt]
\hline
Convergence rate & \textemdash & 0.986914 & 0.972332 & 0.963969 & 0.963169 \\[3pt]
\hline
\end{tabular}}
\label{tab1}
\end{table}

 \begin{exa} \label{example2}{\bf ($100$ dimensional problem without jump)} 
 We now consider the problem without jump process defined on the 100-dimensional hypercube domain $\Omega = (-1,1)^{100}$ to further assess the applicability of the proposed algorithm to high-dimensional settings. In this case, the non-divergence form operator $\mathcal{L}$ in \eqref{mainprob} is given by
 
\begin{equation*}
\mathcal{L}[u](t,\bx)  = \frac{1}{2}{\rm Trace}\big( \sigma(t,\bx) \sigma(t,\bx)^{\top}{\rm Hess}_{\bx}u(t,\bx) \big) + \big\langle  \mu(t,\bx), \nabla_{\bx}u(t,\bx) \big\rangle_{\R^{d}},
\end{equation*}
 where the drift $\mu(t,\bx) =\big(\mu(t,x_{1}),\cdots, \mu(t,x_{d})\big)^{\top},$  $\mu(t,x_{i}) =  \frac{\sin(t)}{1+t^{2}}e^{-x_{i}^{2}}$, and diffusion coefficients are
\begin{eqnarray*}
 \begin{split}
&\sigma(t,\bx) = e^{-5t}\begin{pmatrix}
  \cos(x_{1} + x_{2}) &\cdots&0&0 \\
 \vdots&  \ddots& \vdots&\vdots\\[6pt]
 0&\cdots& \cos(x_{d-1} + x_{d}) &0\\[6pt]
  0&\cdots& 0&\sin(x_{d})
\end{pmatrix}. \end{split}
 \end{eqnarray*}
The problem admits the exact solution
 \begin{equation}\label{exactsolu2}
 u(t,\bx) =  \cos(t^{2})e^{-|\bx|^{2}},\;\;\;(t,\bx)\in [0,\infty)\times\Omega,
 \end{equation}
which determines the initial condition $u^0(\bx)$. An additional source term is added to the equation such that the resulting nonlinear forcing term takes the form $f(u)={\rm e}^{-2|u|}$.  
We set the terminal time $T=1$, and numerically solve \eqref{mainprob} on the cubic domain $(-1,1)^{100}$. 
\end{exa}

\begin{figure}[!h]\hspace{0pt}
\subfigure{
\begin{minipage}[t]{0.4\textwidth}
\centering
\rotatebox[origin=cc]{-0}{\includegraphics[width=1.2\textwidth,height=1\textwidth]{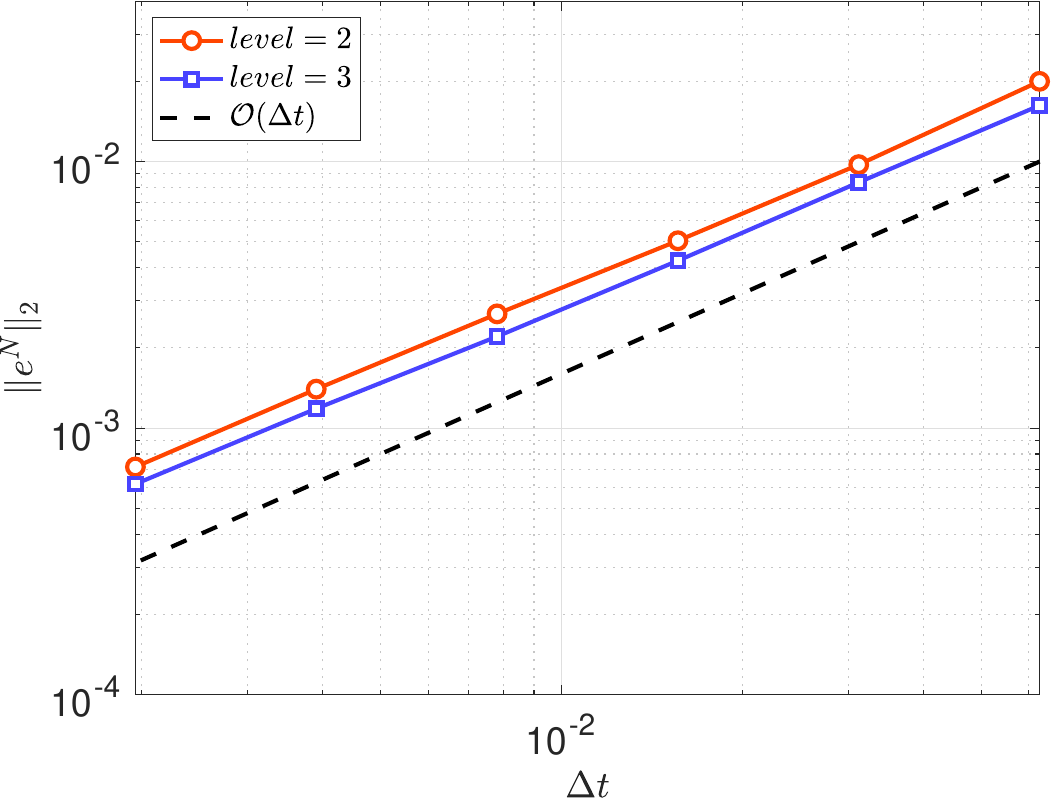}}
\end{minipage}}\hspace{30pt}
\subfigure{
\begin{minipage}[t]{0.4\textwidth}
\centering
\rotatebox[origin=cc]{-0}{\includegraphics[width=1.2\textwidth,height=1\textwidth]{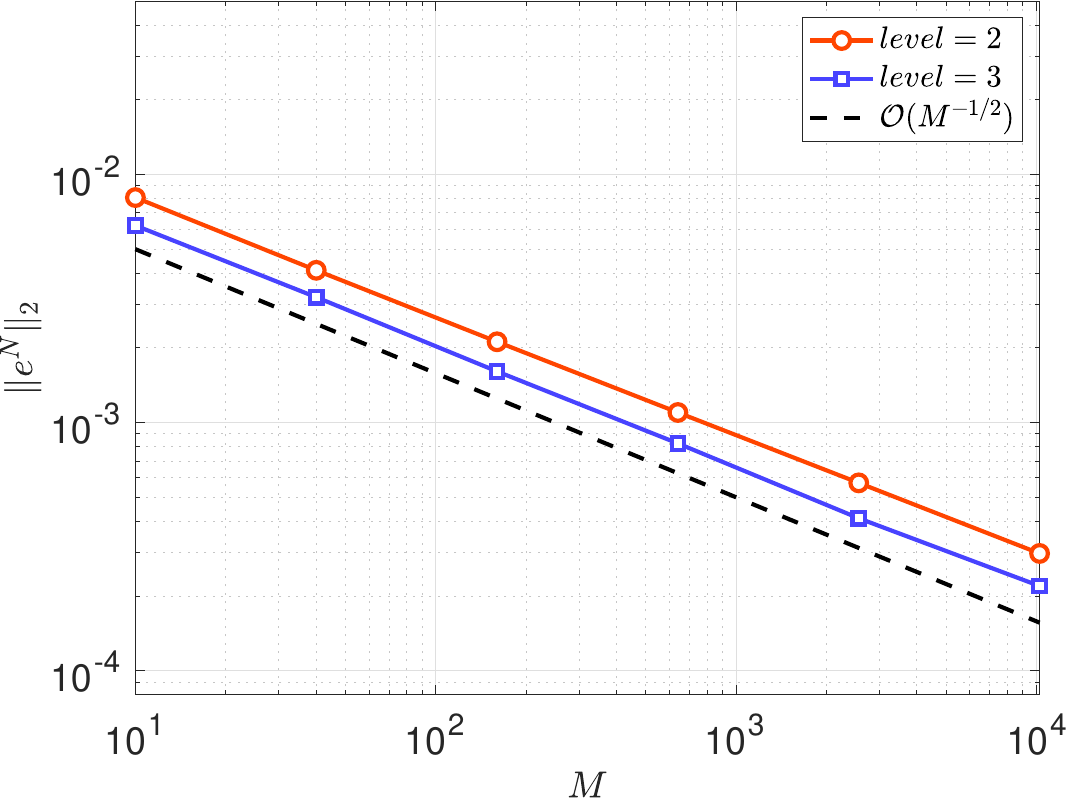}}
\end{minipage}}
\vskip -5pt
\caption
{\small  Left:  The numerical errors of \eqref{exactsolu1} against various time step $\Delta t$ with $T=1$ and $M=10000$; Right: The numerical errors of \eqref{exactsolu1} against the number of the path $M$ with $T=1$ and $\Delta t=2^{-10}$.}\label{err_100D}
\end{figure}

In this example, we perform simulations on two sparse grids with levels $\bm{l} = 2$ and $\bm{l} = 3$ to assess the accuracy and efficiency of the proposed algorithm under different spatial resolutions. Figure~\ref{err_100D} shows the $L^2$ errors for different time step sizes $\Delta t$ and numbers of Monte Carlo samples $M$. As shown in Figure~\ref{err_100D} (left), the error decreases linearly with respect to $\Delta t$, which confirm the first-order convergence in time. This observation is consistent with the theoretical predictions and validates the temporal accuracy of the scheme. Figure~\ref{err_100D} (right) illustrates the $L^2$ error as a function of $M$ for different grid levels. The results demonstrate a convergence rate of approximately $\mathcal{O}(M^{-1/2})$, which matches the expected behavior of Monte Carlo methods. In addition, for fixed $\Delta t$ and $M$, increasing the sparse grid level $\bm{l}$ leads to a clear reduction in the total error. This improvement arises from the enhanced spatial resolution at higher levels, which yields more accurate interpolation and solution reconstruction. These findings underscore the effectiveness of sparse grid techniques in preserving accuracy in high-dimensional problems.

\begin{exa}\label{exa3}{\bf(100 dimensional problem with hypersingular kernel functions)}
Finally, we consider problem \eqref{mainprob} on the 100-dimensional domain $\Omega = (-1, 1)^{100}$ with a hypersingular kernel.  
We take the hypersingular kernel $\varphi(\bz) = C_H |\bz|^{-d-\alpha}$, where the constant $C_H$ is defined in \eqref{kernelsingular}. The nonlinear forcing term is set to $f(u) = |\sin(u)|^{2/3}$,  the nonlocal boundary condition is given by $g(t, \bx) = 0$, $\mu(t,\bx) = \big(\mu(t,x_{1}),\cdots, \mu(t,x_{d})\big)^{\top},$ $\mu(t,x_{i}) =  e^{-t^{2}} \ln\Big(3+\frac{|x_{i}|}{1+t}\Big),$ $c(t,\bx,\bz) =\Big(e^{-z_{1}} + t x_{1},e^{-z_{2}}+\frac12 tx_{2} ,\cdots, e^{-z_{d-1}}+\frac{1}{d-1} tx_{d-1} ,e^{-z_{d}}+\frac1d tx_{d} \Big)^{\top}$, and  
\begin{equation*}
 \begin{split}
&\sigma(t,\bx) = \frac{\cos(t)}{1+10t^{2}}\begin{pmatrix}
 \sin( x_{1} \cdot \frac 1 2 x^{2}_{2}) &\cdots&0&0 \\
  0 &\sin( x_{2} \cdot \frac 1 3 x^{3}_{3}) &0&0\\
 \vdots&  \ddots& \vdots&\vdots\\[6pt]
 0&\cdots& \sin( x_{d-1} \cdot \frac 1 d x^{d}_{d}) &0\\[6pt]
  0&\cdots& 0&\cos(x_{1}\cdot x_{d})
\end{pmatrix}.
 \end{split}
 \end{equation*}
\end{exa}
\begin{figure}[!h]\hspace{0pt}
\subfigure{
\begin{minipage}[t]{0.4\textwidth}
\centering
\rotatebox[origin=cc]{-0}{\includegraphics[width=1.2\textwidth,height=1\textwidth]{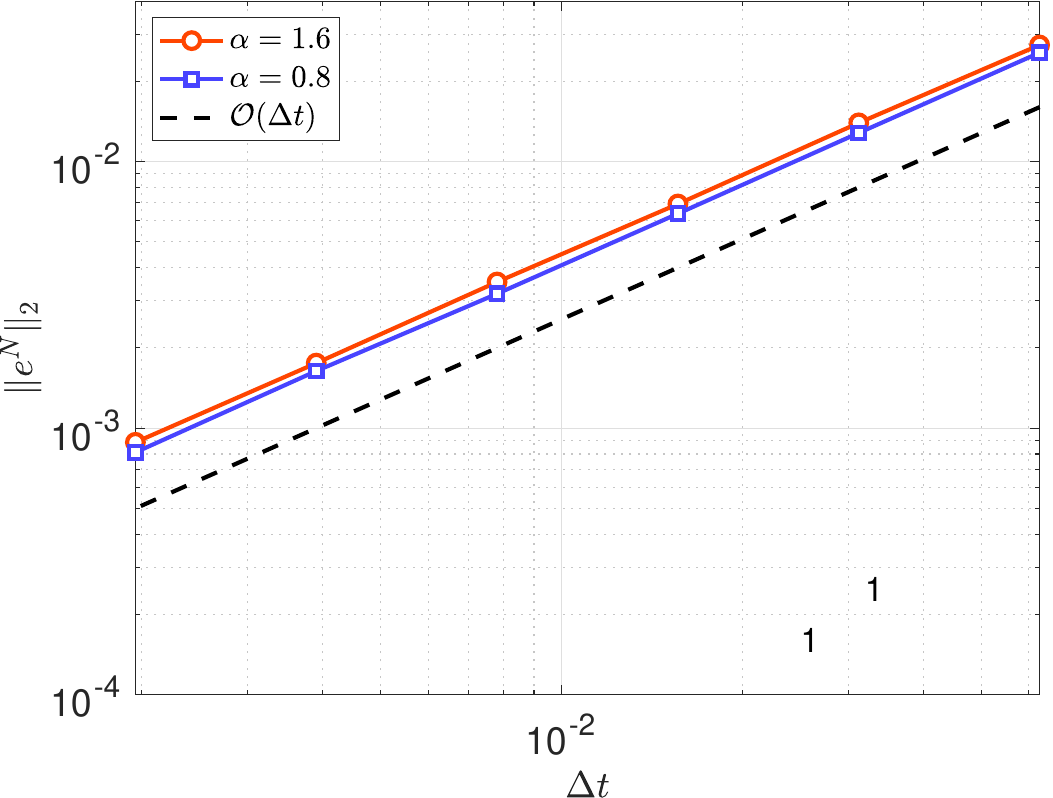}}
\end{minipage}}\hspace{30pt}
\subfigure{
\begin{minipage}[t]{0.4\textwidth}
\centering
\rotatebox[origin=cc]{-0}{\includegraphics[width=1.2\textwidth,height=1\textwidth]{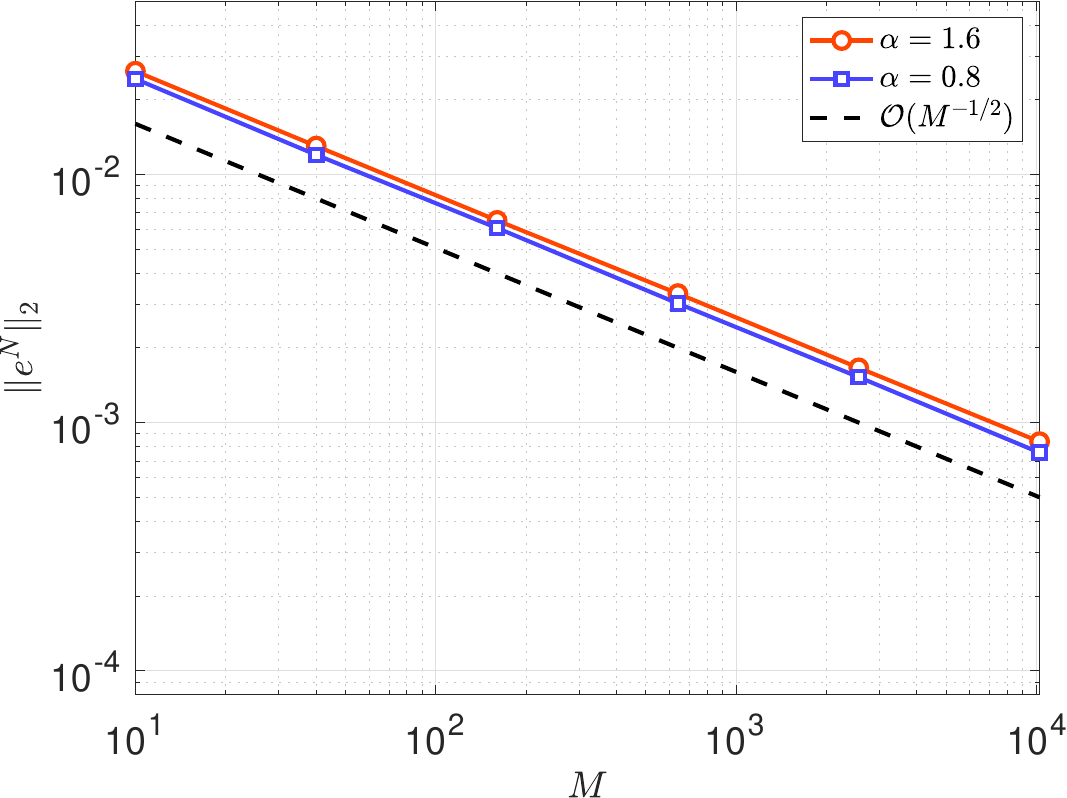}}
\end{minipage}}
\vskip -5pt
\caption
{\small  Left:  The numerical errors of Example \ref{exa3} against various time step $\Delta t$ with $T=1$ and $M=10000$; Right: The numerical errors of Example \ref{exa3}     against the number of the path $M$ with $T=1$ and $\Delta t=2^{-10}$.}\label{err_100D_q}
\end{figure}
In this example, since no analytical solution is available, we construct a reference solution by employing a refined time step $\Delta t = 2^{-11}$ and a sufficiently large number of Monte Carlo paths $M = 100000$. This solution serves as a benchmark for evaluating the convergence rate of the proposed algorithm. A sparse grid with level $\bm{l} = 2$ is used to assess the error behavior under varying time step sizes $\Delta t$ and numbers of simulation paths $M$. The numerical results are presented in Figure~\ref{err_100D_q}. As shown in the figure, the $L^2$-error exhibits first-order convergence with respect to $\Delta t$ and approximately half-order convergence with respect to $M$, consistent with the behavior observed in Example 2. Furthermore, we include results for different values of the fractional parameter $\alpha$ in the hypersingular kernel and find that the convergence rates remain unaffected by the choice of $\alpha$. Notably, we also observe a similar phenomenon reported in \cite{SSX2025}: the numerical error decreases as the fractional order $\alpha$ becomes smaller, underscoring a key advantage of the proposed method over traditional deterministic approaches in the presence of strong nonlocality.

\section{Conclusion}
This paper proposes a sparse grid-based implicit-explicit Monte Carlo method for solving high-dimensional semi-linear nonlocal diffusion equations with volume constraints on bounded domains. The method is developed based on the Feynman-Kac representation and is notable for its structural simplicity and ease of implementation. The key idea lies in simulating stochastic processes to track the evolution of the solution over each time subinterval, while employing an implicit-explicit time-stepping scheme to handle the nonlinearity of the associated PIDEs.
A rigorous error analysis is provided, showing that the accuracy of the method is primarily governed by the time step size and the number of Monte Carlo samples, with convergence rates of order $\mathcal{O}(\Delta t)$ and $\mathcal{O}(M^{-1/2})$, respectively. Extensive numerical experiments are conducted to validate the theoretical findings. The numerical results demonstrate that the proposed algorithm is particularly effective in high-dimensional settings, owing to the synergistic integration of sparse grid interpolation, stochastic simulation, and the implicit-explicit scheme. Moreover, the method is inherently parallelizable, making it highly scalable and well-suited for large-scale computations. In summary, the proposed approach offers a robust and efficient computational framework for high-dimensional semi-linear nonlocal PDEs, combining theoretical analysis with practical efficiency through sparse grid techniques and parallelism, and is promising for a wide range of real-world applications.

 \bigskip
 
\noindent{\bf Acknowledgements} \hspace{4pt} The first author was supported in part by the National Natural Science Foundation of China (Nos. 12201385 and 12271365) and the Fundamental Research Funds for the Central Universities 2021110474. The second author was supported by the Scientific Research Start-up Funds of Hainan University  XJ2400010491.


\bibliography{sn-bibliography.bib}	


\begin{thebibliography}{67}
\ifx \bisbn   \undefined \def \bisbn  #1{ISBN #1}\fi
\ifx \binits  \undefined \def \binits#1{#1}\fi
\ifx \bauthor  \undefined \def \bauthor#1{#1}\fi
\ifx \batitle  \undefined \def \batitle#1{#1}\fi
\ifx \bjtitle  \undefined \def \bjtitle#1{#1}\fi
\ifx \bvolume  \undefined \def \bvolume#1{\textbf{#1}}\fi
\ifx \byear  \undefined \def \byear#1{#1}\fi
\ifx \bissue  \undefined \def \bissue#1{#1}\fi
\ifx \bfpage  \undefined \def \bfpage#1{#1}\fi
\ifx \blpage  \undefined \def \blpage #1{#1}\fi
\ifx \burl  \undefined \def \burl#1{\textsf{#1}}\fi
\ifx \doiurl  \undefined \def \doiurl#1{\url{https://doi.org/#1}}\fi
\ifx \betal  \undefined \def \betal{\textit{et al.}}\fi
\ifx \binstitute  \undefined \def \binstitute#1{#1}\fi
\ifx \binstitutionaled  \undefined \def \binstitutionaled#1{#1}\fi
\ifx \bctitle  \undefined \def \bctitle#1{#1}\fi
\ifx \beditor  \undefined \def \beditor#1{#1}\fi
\ifx \bpublisher  \undefined \def \bpublisher#1{#1}\fi
\ifx \bbtitle  \undefined \def \bbtitle#1{#1}\fi
\ifx \bedition  \undefined \def \bedition#1{#1}\fi
\ifx \bseriesno  \undefined \def \bseriesno#1{#1}\fi
\ifx \blocation  \undefined \def \blocation#1{#1}\fi
\ifx \bsertitle  \undefined \def \bsertitle#1{#1}\fi
\ifx \bsnm \undefined \def \bsnm#1{#1}\fi
\ifx \bsuffix \undefined \def \bsuffix#1{#1}\fi
\ifx \bparticle \undefined \def \bparticle#1{#1}\fi
\ifx \barticle \undefined \def \barticle#1{#1}\fi
\bibcommenthead
\ifx \bconfdate \undefined \def \bconfdate #1{#1}\fi
\ifx \botherref \undefined \def \botherref #1{#1}\fi
\ifx \url \undefined \def \url#1{\textsf{#1}}\fi
\ifx \bchapter \undefined \def \bchapter#1{#1}\fi
\ifx \bbook \undefined \def \bbook#1{#1}\fi
\ifx \bcomment \undefined \def \bcomment#1{#1}\fi
\ifx \oauthor \undefined \def \oauthor#1{#1}\fi
\ifx \citeauthoryear \undefined \def \citeauthoryear#1{#1}\fi
\ifx \endbibitem  \undefined \def \endbibitem {}\fi
\ifx \bconflocation  \undefined \def \bconflocation#1{#1}\fi
\ifx \arxivurl  \undefined \def \arxivurl#1{\textsf{#1}}\fi
\csname PreBibitemsHook\endcsname

\bibitem[\protect\citeauthoryear{Del-Castillo-Negrete}{2010}]{del2010}
\begin{barticle}
\bauthor{\bsnm{Del-Castillo-Negrete}, \binits{D.}}:
\batitle{Non-diffusive, non-local transport in fluids and plasmas}.
\bjtitle{Nonlinear Processes in Geophysics}
\bvolume{17}(\bissue{6}),
\bfpage{795}--\blpage{807}
(\byear{2010})
\end{barticle}
\endbibitem

\bibitem[\protect\citeauthoryear{Del-Castillo-Negrete}{2006}]{del2006}
\begin{botherref}
\oauthor{\bsnm{Del-Castillo-Negrete}, \binits{D.}}:
Fractional diffusion models of nonlocal transport.
Physics of Plasmas
\textbf{13}(8)
(2006)
\end{botherref}
\endbibitem

\bibitem[\protect\citeauthoryear{S{\'a}nchez et~al.}{2008}]{Sanchez2008}
\begin{barticle}
\bauthor{\bsnm{S{\'a}nchez}, \binits{R.}},
\bauthor{\bsnm{Newman}, \binits{D.E.}},
\bauthor{\bsnm{Leboeuf}, \binits{J.-N.}},
\bauthor{\bsnm{Decyk}, \binits{V.}},
\bauthor{\bsnm{Carreras}, \binits{B.A.}}:
\batitle{Nature of transport across sheared zonal flows in electrostatic
  ion-temperature-gradient gyrokinetic plasma turbulence}.
\bjtitle{Physical Review Letters}
\bvolume{101}(\bissue{20}),
\bfpage{205002}
(\byear{2008})
\end{barticle}
\endbibitem

\bibitem[\protect\citeauthoryear{Chen}{2006}]{chen2006}
\begin{botherref}
\oauthor{\bsnm{Chen}, \binits{W.}}:
A speculative study of 2/ 3-order fractional {L}aplacian modeling of
  turbulence: {S}ome thoughts and conjectures.
Chaos: An Interdisciplinary Journal of Nonlinear Science
\textbf{16}(2)
(2006)
\end{botherref}
\endbibitem

\bibitem[\protect\citeauthoryear{Gunzburger et~al.}{2018}]{Gunzburger2018}
\begin{barticle}
\bauthor{\bsnm{Gunzburger}, \binits{M.}},
\bauthor{\bsnm{Jiang}, \binits{N.}},
\bauthor{\bsnm{Xu}, \binits{F.}}:
\batitle{Analysis and approximation of a fractional {L}aplacian-based closure
  model for turbulent flows and its connection to richardson pair dispersion}.
\bjtitle{Computers \& Mathematics with Applications}
\bvolume{75}(\bissue{6}),
\bfpage{1973}--\blpage{2001}
(\byear{2018})
\end{barticle}
\endbibitem

\bibitem[\protect\citeauthoryear{Silling}{2000}]{Silling2000}
\begin{barticle}
\bauthor{\bsnm{Silling}, \binits{S.A.}}:
\batitle{Reformulation of elasticity theory for discontinuities and long-range
  forces}.
\bjtitle{Journal of the Mechanics and Physics of Solids}
\bvolume{48}(\bissue{1}),
\bfpage{175}--\blpage{209}
(\byear{2000})
\end{barticle}
\endbibitem

\bibitem[\protect\citeauthoryear{Silling and Askari}{2005}]{Silling2005}
\begin{barticle}
\bauthor{\bsnm{Silling}, \binits{S.A.}},
\bauthor{\bsnm{Askari}, \binits{E.}}:
\batitle{A meshfree method based on the peridynamic model of solid mechanics}.
\bjtitle{Computers \& Structures}
\bvolume{83}(\bissue{17-18}),
\bfpage{1526}--\blpage{1535}
(\byear{2005})
\end{barticle}
\endbibitem

\bibitem[\protect\citeauthoryear{Tao and Winkler}{2011}]{tao2011}
\begin{barticle}
\bauthor{\bsnm{Tao}, \binits{Y.}},
\bauthor{\bsnm{Winkler}, \binits{M.}}:
\batitle{A chemotaxis-haptotaxis model: the roles of nonlinear diffusion and
  logistic source}.
\bjtitle{SIAM Journal on Mathematical Analysis}
\bvolume{43}(\bissue{2}),
\bfpage{685}--\blpage{704}
(\byear{2011})
\end{barticle}
\endbibitem

\bibitem[\protect\citeauthoryear{Ahmed and Elgazzar}{2007}]{Ahmed2007}
\begin{barticle}
\bauthor{\bsnm{Ahmed}, \binits{E.}},
\bauthor{\bsnm{Elgazzar}, \binits{A.}}:
\batitle{On fractional order differential equations model for nonlocal
  epidemics}.
\bjtitle{Physica A: Statistical Mechanics and its Applications}
\bvolume{379}(\bissue{2}),
\bfpage{607}--\blpage{614}
(\byear{2007})
\end{barticle}
\endbibitem

\bibitem[\protect\citeauthoryear{Gilboa and Osher}{2009}]{Gilboa2009}
\begin{barticle}
\bauthor{\bsnm{Gilboa}, \binits{G.}},
\bauthor{\bsnm{Osher}, \binits{S.}}:
\batitle{Nonlocal operators with applications to image processing}.
\bjtitle{Multiscale Modeling \& Simulation}
\bvolume{7}(\bissue{3}),
\bfpage{1005}--\blpage{1028}
(\byear{2009})
\end{barticle}
\endbibitem

\bibitem[\protect\citeauthoryear{Buades et~al.}{2010}]{Buades2010}
\begin{barticle}
\bauthor{\bsnm{Buades}, \binits{A.}},
\bauthor{\bsnm{Coll}, \binits{B.}},
\bauthor{\bsnm{Morel}, \binits{J.-M.}}:
\batitle{Image denoising methods. {A} new nonlocal principle}.
\bjtitle{SIAM Review}
\bvolume{52}(\bissue{1}),
\bfpage{113}--\blpage{147}
(\byear{2010})
\end{barticle}
\endbibitem

\bibitem[\protect\citeauthoryear{Cont and Voltchkova}{2005}]{cont2005integro}
\begin{barticle}
\bauthor{\bsnm{Cont}, \binits{R.}},
\bauthor{\bsnm{Voltchkova}, \binits{E.}}:
\batitle{Integro-differential equations for option prices in exponential
  {L}{\'e}vy models}.
\bjtitle{Finance and Stochastics}
\bvolume{9},
\bfpage{299}--\blpage{325}
(\byear{2005})
\end{barticle}
\endbibitem

\bibitem[\protect\citeauthoryear{Duo and Zhang}{2019}]{duo2018finite}
\begin{barticle}
\bauthor{\bsnm{Duo}, \binits{S.}},
\bauthor{\bsnm{Zhang}, \binits{Y.}}:
\batitle{Accurate numerical methods for two and three dimensional integral
  fractional {L}aplacian with applications}.
\bjtitle{Computer Methods in Applied Mechanics and Engineering}
\bvolume{355},
\bfpage{639}--\blpage{662}
(\byear{2019})
\end{barticle}
\endbibitem

\bibitem[\protect\citeauthoryear{Hao et~al.}{2021}]{hao2019fractional}
\begin{barticle}
\bauthor{\bsnm{Hao}, \binits{Z.}},
\bauthor{\bsnm{Zhang}, \binits{Z.}},
\bauthor{\bsnm{Du}, \binits{R.}}:
\batitle{Fractional centered difference scheme for high-dimensional integral
  fractional {L}aplacian}.
\bjtitle{Journal of Computational Physics}
\bvolume{424},
\bfpage{109851}
(\byear{2021})
\end{barticle}
\endbibitem

\bibitem[\protect\citeauthoryear{Minden and Ying}{2020}]{victor2020simple}
\begin{barticle}
\bauthor{\bsnm{Minden}, \binits{V.}},
\bauthor{\bsnm{Ying}, \binits{L.}}:
\batitle{A simple solver for the fractional {L}aplacian in multiple
  dimensions}.
\bjtitle{SIAM Journal on Scientific Computing}
\bvolume{42}(\bissue{2}),
\bfpage{878}--\blpage{900}
(\byear{2020})
\end{barticle}
\endbibitem

\bibitem[\protect\citeauthoryear{Acosta and
  Borthagaray}{2017}]{acosta2017fractional}
\begin{barticle}
\bauthor{\bsnm{Acosta}, \binits{G.}},
\bauthor{\bsnm{Borthagaray}, \binits{J.P.}}:
\batitle{A fractional {L}aplace equation: regularity of solutions and finite
  element approximations}.
\bjtitle{SIAM Journal on Numerical Analysis}
\bvolume{55}(\bissue{2}),
\bfpage{472}--\blpage{495}
(\byear{2017})
\end{barticle}
\endbibitem

\bibitem[\protect\citeauthoryear{Ainsworth and Glusa}{2017}]{Ainsworth2017}
\begin{barticle}
\bauthor{\bsnm{Ainsworth}, \binits{M.}},
\bauthor{\bsnm{Glusa}, \binits{C.}}:
\batitle{Aspects of an adaptive finite element method for the fractional
  laplacian: a priori and a posteriori error estimates, efficient
  implementation and multigrid solver}.
\bjtitle{Computer Methods in Applied Mechanics and Engineering}
\bvolume{327},
\bfpage{4}--\blpage{35}
(\byear{2017})
\end{barticle}
\endbibitem

\bibitem[\protect\citeauthoryear{Bonito et~al.}{2019}]{bonito2019numerical}
\begin{barticle}
\bauthor{\bsnm{Bonito}, \binits{A.}},
\bauthor{\bsnm{Lei}, \binits{W.}},
\bauthor{\bsnm{Pasciak}, \binits{J.E.}}:
\batitle{Numerical approximation of the integral fractional {L}aplacian}.
\bjtitle{Numerische Mathematik}
\bvolume{142}(\bissue{2}),
\bfpage{235}--\blpage{278}
(\byear{2019})
\end{barticle}
\endbibitem

\bibitem[\protect\citeauthoryear{Mao and Shen}{2017}]{mao2017hermite}
\begin{barticle}
\bauthor{\bsnm{Mao}, \binits{Z.}},
\bauthor{\bsnm{Shen}, \binits{J.}}:
\batitle{Hermite spectral methods for fractional {PDE}s in unbounded domains}.
\bjtitle{SIAM Journal on Scientific Computing}
\bvolume{39}(\bissue{5}),
\bfpage{1928}--\blpage{1950}
(\byear{2017})
\end{barticle}
\endbibitem

\bibitem[\protect\citeauthoryear{Tang et~al.}{2018}]{tang2018hermite}
\begin{botherref}
\oauthor{\bsnm{Tang}, \binits{T.}},
\oauthor{\bsnm{Yuan}, \binits{H.}},
\oauthor{\bsnm{Zhou}, \binits{T.}}:
Hermite spectral collocation methods for fractional {PDE}s in unbounded
  domains.
arXiv preprint arXiv:1801.09073
(2018)
\end{botherref}
\endbibitem

\bibitem[\protect\citeauthoryear{Sheng et~al.}{2020}]{sheng2019fast}
\begin{barticle}
\bauthor{\bsnm{Sheng}, \binits{C.}},
\bauthor{\bsnm{Shen}, \binits{J.}},
\bauthor{\bsnm{Tang}, \binits{T.}},
\bauthor{\bsnm{Wang}, \binits{L.-L.}},
\bauthor{\bsnm{Yuan}, \binits{H.}}:
\batitle{Fast {F}ourier-like mapped {C}hebyshev spectral-{G}alerkin methods for
  {PDE}s with integral fractional {L}aplacian in unbounded domains}.
\bjtitle{SIAM Journal on Numerical Analysis}
\bvolume{58}(\bissue{5}),
\bfpage{2435}--\blpage{2464}
(\byear{2020})
\end{barticle}
\endbibitem

\bibitem[\protect\citeauthoryear{Tian and Du}{2013}]{Tian2013}
\begin{barticle}
\bauthor{\bsnm{Tian}, \binits{X.}},
\bauthor{\bsnm{Du}, \binits{Q.}}:
\batitle{Analysis and comparison of different approximations to nonlocal
  diffusion and linear peridynamic equations}.
\bjtitle{SIAM Journal on Numerical Analysis}
\bvolume{51}(\bissue{6}),
\bfpage{3458}--\blpage{3482}
(\byear{2013})
\end{barticle}
\endbibitem

\bibitem[\protect\citeauthoryear{Leng et~al.}{2021}]{Leng2021}
\begin{barticle}
\bauthor{\bsnm{Leng}, \binits{Y.}},
\bauthor{\bsnm{Tian}, \binits{X.}},
\bauthor{\bsnm{Trask}, \binits{N.}},
\bauthor{\bsnm{Foster}, \binits{J.T.}}:
\batitle{Asymptotically compatible reproducing kernel collocation and meshfree
  integration for nonlocal diffusion}.
\bjtitle{SIAM Journal on Numerical Analysis}
\bvolume{59}(\bissue{1}),
\bfpage{88}--\blpage{118}
(\byear{2021})
\end{barticle}
\endbibitem

\bibitem[\protect\citeauthoryear{Aggarwal et~al.}{2024}]{Aggarwal2024}
\begin{barticle}
\bauthor{\bsnm{Aggarwal}, \binits{A.}},
\bauthor{\bsnm{Holden}, \binits{H.}},
\bauthor{\bsnm{Vaidya}, \binits{G.}}:
\batitle{On the accuracy of the finite volume approximations to nonlocal
  conservation laws}.
\bjtitle{Numerische Mathematik}
\bvolume{156}(\bissue{1}),
\bfpage{237}--\blpage{271}
(\byear{2024})
\end{barticle}
\endbibitem

\bibitem[\protect\citeauthoryear{D’Elia et~al.}{2020}]{Elia2020}
\begin{barticle}
\bauthor{\bsnm{D’Elia}, \binits{M.}},
\bauthor{\bsnm{Du}, \binits{Q.}},
\bauthor{\bsnm{Glusa}, \binits{C.}},
\bauthor{\bsnm{Gunzburger}, \binits{M.}},
\bauthor{\bsnm{Tian}, \binits{X.}},
\bauthor{\bsnm{Zhou}, \binits{Z.}}:
\batitle{Numerical methods for nonlocal and fractional models}.
\bjtitle{Acta Numerica}
\bvolume{29},
\bfpage{1}--\blpage{124}
(\byear{2020})
\end{barticle}
\endbibitem

\bibitem[\protect\citeauthoryear{Du et~al.}{2023}]{Du2023}
\begin{bchapter}
\bauthor{\bsnm{Du}, \binits{Q.}},
\bauthor{\bsnm{Tian}, \binits{X.}},
\bauthor{\bsnm{Zhou}, \binits{Z.}}:
\bctitle{Nonlocal diffusion models with consistent local and fractional
  limits}.
In: \bbtitle{A$^3$N$^2$M: Approximation, Applications, and Analysis of
  Nonlocal, Nonlinear Models: Proceedings of the 50th John H. Barrett Memorial
  Lectures},
pp. \bfpage{175}--\blpage{213}
(\byear{2023}).
\bcomment{Springer}
\end{bchapter}
\endbibitem

\bibitem[\protect\citeauthoryear{Du}{2019}]{Du2019}
\begin{bbook}
\bauthor{\bsnm{Du}, \binits{Q.}}:
\bbtitle{Nonlocal {M}odeling, {A}nalysis, and {C}omputation: {N}onlocal
  {M}odeling, {A}nalysis, And {C}omputation},
(\byear{2019}).
\bcomment{SIAM}
\end{bbook}
\endbibitem

\bibitem[\protect\citeauthoryear{Pang et~al.}{2020}]{PANG2020109760}
\begin{barticle}
\bauthor{\bsnm{Pang}, \binits{G.}},
\bauthor{\bsnm{D'Elia}, \binits{M.}},
\bauthor{\bsnm{Parks}, \binits{M.}},
\bauthor{\bsnm{Karniadakis}, \binits{G.E.}}:
\batitle{npinns: {N}onlocal physics-informed neural networks for a parametrized
  nonlocal universal laplacian operator. {A}lgorithms and applications}.
\bjtitle{Journal of Computational Physics}
\bvolume{422},
\bfpage{109760}
(\byear{2020})
\end{barticle}
\endbibitem

\bibitem[\protect\citeauthoryear{Guo et~al.}{2022}]{Guo2022M}
\begin{barticle}
\bauthor{\bsnm{Guo}, \binits{L.}},
\bauthor{\bsnm{Wu}, \binits{H.}},
\bauthor{\bsnm{Yu}, \binits{X.}},
\bauthor{\bsnm{Zhou}, \binits{T.}}:
\batitle{Monte {C}arlo fpinns: {D}eep learning method for forward and inverse
  problems involving high dimensional fractional partial differential
  equations}.
\bjtitle{Computer Methods in Applied Mechanics and Engineering}
\bvolume{400},
\bfpage{115523}
(\byear{2022})
\end{barticle}
\endbibitem

\bibitem[\protect\citeauthoryear{Castro}{2022}]{Castro2022}
\begin{barticle}
\bauthor{\bsnm{Castro}, \binits{J.}}:
\batitle{Deep learning schemes for parabolic nonlocal integro-differential
  equations}.
\bjtitle{Partial Differ. Equ. Appl.}
\bvolume{3}(\bissue{6}),
\bfpage{77}--\blpage{35}
(\byear{2022})
\end{barticle}
\endbibitem

\bibitem[\protect\citeauthoryear{Lu et~al.}{2024}]{Lu2024}
\begin{barticle}
\bauthor{\bsnm{Lu}, \binits{L.}},
\bauthor{\bsnm{Guo}, \binits{H.}},
\bauthor{\bsnm{Yang}, \binits{X.}},
\bauthor{\bsnm{Zhu}, \binits{Y.}}:
\batitle{Temporal difference learning for high-dimensional {PIDE}s with jumps}.
\bjtitle{SIAM Journal on Scientific Computing}
\bvolume{46}(\bissue{4}),
\bfpage{349}--\blpage{368}
(\byear{2024})
\end{barticle}
\endbibitem

\bibitem[\protect\citeauthoryear{Du et~al.}{2014}]{DuHuang2014}
\begin{botherref}
\oauthor{\bsnm{Du}, \binits{Q.}},
\oauthor{\bsnm{Huang}, \binits{Z.}},
\oauthor{\bsnm{Lehoucq}, \binits{R.B.}}:
Nonlocal convection-diffusion volume-constrained problems and jump processes.
Discrete \& Continuous Dynamical Systems-Series B
\textbf{19}(2)
(2014)
\end{botherref}
\endbibitem

\bibitem[\protect\citeauthoryear{Metzler and Klafter}{2000}]{Metzler2000}
\begin{barticle}
\bauthor{\bsnm{Metzler}, \binits{R.}},
\bauthor{\bsnm{Klafter}, \binits{J.}}:
\batitle{The random walk's guide to anomalous diffusion: a fractional dynamics
  approach}.
\bjtitle{Physics Reports}
\bvolume{339}(\bissue{1}),
\bfpage{1}--\blpage{77}
(\byear{2000})
\end{barticle}
\endbibitem

\bibitem[\protect\citeauthoryear{Pardoux and Peng}{2005}]{Pardoux2005}
\begin{bchapter}
\bauthor{\bsnm{Pardoux}, \binits{E.}},
\bauthor{\bsnm{Peng}, \binits{S.}}:
\bctitle{Backward stochastic differential equations and quasilinear parabolic
  partial differential equations}.
In: \bbtitle{Stochastic Partial Differential Equations and Their Applications:
  Proceedings of IFIP WG 7/1 International Conference University of North
  Carolina at Charlotte, NC June 6--8, 1991},
pp. \bfpage{200}--\blpage{217}
(\byear{2005}).
\bcomment{Springer}
\end{bchapter}
\endbibitem

\bibitem[\protect\citeauthoryear{Robbe et~al.}{2017}]{Robbe2017NM}
\begin{barticle}
\bauthor{\bsnm{Robbe}, \binits{P.}},
\bauthor{\bsnm{Nuyens}, \binits{D.}},
\bauthor{\bsnm{Vandewalle}, \binits{S.}}:
\batitle{A multi-index quasi--{M}onte {C}arlo algorithm for lognormal diffusion
  problems}.
\bjtitle{SIAM Journal on Scientific Computing}
\bvolume{39}(\bissue{5}),
\bfpage{851}--\blpage{872}
(\byear{2017})
\end{barticle}
\endbibitem

\bibitem[\protect\citeauthoryear{L{\o}vbak et~al.}{2021}]{Lovbak2021NM}
\begin{barticle}
\bauthor{\bsnm{L{\o}vbak}, \binits{E.}},
\bauthor{\bsnm{Samaey}, \binits{G.}},
\bauthor{\bsnm{Vandewalle}, \binits{S.}}:
\batitle{A multilevel {M}onte {C}arlo method for asymptotic-preserving particle
  schemes in the diffusive limit}.
\bjtitle{Numerische Mathematik}
\bvolume{148}(\bissue{1}),
\bfpage{141}--\blpage{186}
(\byear{2021})
\end{barticle}
\endbibitem

\bibitem[\protect\citeauthoryear{Shao and Xiong}{2020}]{Shao2020}
\begin{barticle}
\bauthor{\bsnm{Shao}, \binits{S.}},
\bauthor{\bsnm{Xiong}, \binits{Y.}}:
\batitle{Branching random walk solutions to the {W}igner equation}.
\bjtitle{SIAM Journal on Numerical Analysis}
\bvolume{58}(\bissue{5}),
\bfpage{2589}--\blpage{2608}
(\byear{2020})
\end{barticle}
\endbibitem

\bibitem[\protect\citeauthoryear{Ding et~al.}{2023}]{Ding2023}
\begin{barticle}
\bauthor{\bsnm{Ding}, \binits{C.}},
\bauthor{\bsnm{Zhou}, \binits{Y.}},
\bauthor{\bsnm{Cai}, \binits{W.}},
\bauthor{\bsnm{Zeng}, \binits{X.}},
\bauthor{\bsnm{Yan}, \binits{C.}}:
\batitle{A path integral {M}onte {C}arlo (pimc) method based on {F}eynman-{K}ac
  formula for electrical impedance tomography}.
\bjtitle{Journal of Computational Physics}
\bvolume{476},
\bfpage{111862}
(\byear{2023})
\end{barticle}
\endbibitem

\bibitem[\protect\citeauthoryear{Lei et~al.}{2025}]{Lei2025an}
\begin{botherref}
\oauthor{\bsnm{Lei}, \binits{Z.}},
\oauthor{\bsnm{Shao}, \binits{S.}},
\oauthor{\bsnm{Xiong}, \binits{Y.}}:
An efficient stochastic particle method for moderately high-dimensional
  nonlinear {PDE}s.
Journal of Computational Physics,
113818
(2025)
\end{botherref}
\endbibitem

\bibitem[\protect\citeauthoryear{Yang et~al.}{2023}]{yang2023}
\begin{barticle}
\bauthor{\bsnm{Yang}, \binits{M.}},
\bauthor{\bsnm{Zhang}, \binits{G.}},
\bauthor{\bsnm{Del-Castillo-Negrete}, \binits{D.}},
\bauthor{\bsnm{Cao}, \binits{Y.}}:
\batitle{A probabilistic scheme for semilinear nonlocal diffusion equations
  with volume constraints}.
\bjtitle{SIAM Journal on Numerical Analysis}
\bvolume{61}(\bissue{6}),
\bfpage{2718}--\blpage{2743}
(\byear{2023})
\end{barticle}
\endbibitem

\bibitem[\protect\citeauthoryear{Kyprianou et~al.}{2018}]{Kyprianou2018}
\begin{barticle}
\bauthor{\bsnm{Kyprianou}, \binits{A.E.}},
\bauthor{\bsnm{Osojnik}, \binits{A.}},
\bauthor{\bsnm{Shardlow}, \binits{T.}}:
\batitle{Unbiased ‘walk-on-spheres’ {M}onte {C}arlo methods for the
  fractional {L}aplacian}.
\bjtitle{IMA Journal of Numerical Analysis}
\bvolume{38}(\bissue{3}),
\bfpage{1550}--\blpage{1578}
(\byear{2018})
\end{barticle}
\endbibitem

\bibitem[\protect\citeauthoryear{Sheng et~al.}{2023}]{SSX2023}
\begin{barticle}
\bauthor{\bsnm{Sheng}, \binits{C.}},
\bauthor{\bsnm{Su}, \binits{B.}},
\bauthor{\bsnm{Xu}, \binits{C.}}:
\batitle{Efficient {M}onte {C}arlo method for integral fractional {L}aplacian
  in multiple dimensions}.
\bjtitle{SIAM Journal on Numerical Analysis}
\bvolume{61}(\bissue{5}),
\bfpage{2035}--\blpage{2061}
(\byear{2023})
\end{barticle}
\endbibitem

\bibitem[\protect\citeauthoryear{Su et~al.}{2023}]{SXS2023}
\begin{barticle}
\bauthor{\bsnm{Su}, \binits{B.}},
\bauthor{\bsnm{Xu}, \binits{C.}},
\bauthor{\bsnm{Sheng}, \binits{C.}}:
\batitle{A new ‘walk on spheres’ type method for fractional diffusion
  equation in high dimensions based on the {F}eynman--{K}ac formulas}.
\bjtitle{Applied Mathematics Letters}
\bvolume{141},
\bfpage{108597}
(\byear{2023})
\end{barticle}
\endbibitem

\bibitem[\protect\citeauthoryear{Tian et~al.}{2016}]{Tian2016}
\begin{barticle}
\bauthor{\bsnm{Tian}, \binits{X.}},
\bauthor{\bsnm{Du}, \binits{Q.}},
\bauthor{\bsnm{Gunzburger}, \binits{M.}}:
\batitle{Asymptotically compatible schemes for the approximation of fractional
  {L}aplacian and related nonlocal diffusion problems on bounded domains}.
\bjtitle{Advances in Computational Mathematics}
\bvolume{42},
\bfpage{1363}--\blpage{1380}
(\byear{2016})
\end{barticle}
\endbibitem

\bibitem[\protect\citeauthoryear{Du et~al.}{2012}]{DuGunzburger2012}
\begin{barticle}
\bauthor{\bsnm{Du}, \binits{Q.}},
\bauthor{\bsnm{Gunzburger}, \binits{M.}},
\bauthor{\bsnm{Lehoucq}, \binits{R.B.}},
\bauthor{\bsnm{Zhou}, \binits{K.}}:
\batitle{Analysis and approximation of nonlocal diffusion problems with volume
  constraints}.
\bjtitle{SIAM Review}
\bvolume{54}(\bissue{4}),
\bfpage{667}--\blpage{696}
(\byear{2012})
\end{barticle}
\endbibitem

\bibitem[\protect\citeauthoryear{Zhang}{2012}]{Zhang2012}
\begin{botherref}
\oauthor{\bsnm{Zhang}, \binits{X.}}:
Stochastic functional differential equations driven by {L}{\'e}vy processes and
  quasi-linear partial integro-differential equations
(2012)
\end{botherref}
\endbibitem

\bibitem[\protect\citeauthoryear{Higham and Kloeden}{2005}]{Higham2005}
\begin{barticle}
\bauthor{\bsnm{Higham}, \binits{D.J.}},
\bauthor{\bsnm{Kloeden}, \binits{P.E.}}:
\batitle{Numerical methods for nonlinear stochastic differential equations with
  jumps}.
\bjtitle{Numerische Mathematik}
\bvolume{101}(\bissue{1}),
\bfpage{101}--\blpage{119}
(\byear{2005})
\end{barticle}
\endbibitem

\bibitem[\protect\citeauthoryear{Barles et~al.}{1997}]{barles1997backward}
\begin{barticle}
\bauthor{\bsnm{Barles}, \binits{G.}},
\bauthor{\bsnm{Buckdahn}, \binits{R.}},
\bauthor{\bsnm{Pardoux}, \binits{E.}}:
\batitle{Backward stochastic differential equations and integral-partial
  differential equations}.
\bjtitle{Stochastics: An International Journal of Probability and Stochastic
  Processes}
\bvolume{60}(\bissue{1-2}),
\bfpage{57}--\blpage{83}
(\byear{1997})
\end{barticle}
\endbibitem

\bibitem[\protect\citeauthoryear{Pardoux and Peng}{1990}]{pardoux1990adapted}
\begin{barticle}
\bauthor{\bsnm{Pardoux}, \binits{E.}},
\bauthor{\bsnm{Peng}, \binits{S.}}:
\batitle{Adapted solution of a backward stochastic differential equation}.
\bjtitle{Systems \& control letters}
\bvolume{14}(\bissue{1}),
\bfpage{55}--\blpage{61}
(\byear{1990})
\end{barticle}
\endbibitem

\bibitem[\protect\citeauthoryear{Kaliuzhnyi-Verbovetskyi and
  Medvedev}{2022}]{Kaliuzhny2022}
\begin{barticle}
\bauthor{\bsnm{Kaliuzhnyi-Verbovetskyi}, \binits{D.}},
\bauthor{\bsnm{Medvedev}, \binits{G.S.}}:
\batitle{Sparse {M}onte {C}arlo method for nonlocal diffusion problems}.
\bjtitle{SIAM Journal on Numerical Analysis}
\bvolume{60}(\bissue{6}),
\bfpage{3001}--\blpage{3028}
(\byear{2022})
\end{barticle}
\endbibitem

\bibitem[\protect\citeauthoryear{Eriksson-Bique et~al.}{2011}]{Eriksson2011NM}
\begin{barticle}
\bauthor{\bsnm{Eriksson-Bique}, \binits{S.}},
\bauthor{\bsnm{Solbrig}, \binits{M.}},
\bauthor{\bsnm{Stefanelli}, \binits{M.}},
\bauthor{\bsnm{Warkentin}, \binits{S.}},
\bauthor{\bsnm{Abbey}, \binits{R.}},
\bauthor{\bsnm{Ipsen}, \binits{I.C.}}:
\batitle{Importance sampling for a {M}onte {C}arlo matrix multiplication
  algorithm, with application to information retrieval}.
\bjtitle{SIAM Journal on Scientific Computing}
\bvolume{33}(\bissue{4}),
\bfpage{1689}--\blpage{1706}
(\byear{2011})
\end{barticle}
\endbibitem

\bibitem[\protect\citeauthoryear{Bungartz and
  Griebel}{2004}]{Bungartz2004NChief}
\begin{barticle}
\bauthor{\bsnm{Bungartz}, \binits{H.-J.}},
\bauthor{\bsnm{Griebel}, \binits{M.}}:
\batitle{Sparse grids}.
\bjtitle{Acta Numerica}
\bvolume{13},
\bfpage{147}--\blpage{269}
(\byear{2004})
\end{barticle}
\endbibitem

\bibitem[\protect\citeauthoryear{Griebel et~al.}{1999}]{Griebel1999NChief}
\begin{barticle}
\bauthor{\bsnm{Griebel}, \binits{M.}},
\bauthor{\bsnm{Oswald}, \binits{P.}},
\bauthor{\bsnm{Schiekofer}, \binits{T.}}:
\batitle{Sparse grids for boundary integral equations}.
\bjtitle{Numerische Mathematik}
\bvolume{83}(\bissue{2}),
\bfpage{279}--\blpage{312}
(\byear{1999})
\end{barticle}
\endbibitem

\bibitem[\protect\citeauthoryear{Smolyak}{1963}]{Smolyak1963}
\begin{bchapter}
\bauthor{\bsnm{Smolyak}, \binits{S.}}:
\bctitle{Quadrature and interpolation formulas for tensor products of certain
  classes of functions}.
In: \bbtitle{Soviet Math. Dokl.},
vol. \bseriesno{4},
pp. \bfpage{240}--\blpage{243}
(\byear{1963})
\end{bchapter}
\endbibitem

\bibitem[\protect\citeauthoryear{Yserentant}{1986}]{Yserentant1986}
\begin{barticle}
\bauthor{\bsnm{Yserentant}, \binits{H.}}:
\batitle{On the multi-level splitting of finite element spaces}.
\bjtitle{Numerische Mathematik}
\bvolume{49}(\bissue{4}),
\bfpage{379}--\blpage{412}
(\byear{1986})
\end{barticle}
\endbibitem

\bibitem[\protect\citeauthoryear{Yserentant}{1987}]{Yserentant1992}
\begin{botherref}
\oauthor{\bsnm{Yserentant}, \binits{H.}}:
Hierarchical bases in the numerical solution of parabolic problems.
Large Scale Scientific Computing,
22--36
(1987)
\end{botherref}
\endbibitem

\bibitem[\protect\citeauthoryear{Temlyakov}{1986}]{Temlyakov1986}
\begin{barticle}
\bauthor{\bsnm{Temlyakov}, \binits{V.N.}}:
\batitle{Approximations of functions with bounded mixed derivative}.
\bjtitle{Trudy Matematicheskogo Instituta imeni VA Steklova}
\bvolume{178},
\bfpage{3}--\blpage{113}
(\byear{1986})
\end{barticle}
\endbibitem

\bibitem[\protect\citeauthoryear{Shen and Yu}{2010}]{SY2010}
\begin{barticle}
\bauthor{\bsnm{Shen}, \binits{J.}},
\bauthor{\bsnm{Yu}, \binits{H.}}:
\batitle{Efficient spectral sparse grid methods and applications to
  high-dimensional elliptic problems}.
\bjtitle{SIAM Journal on Scientific Computing}
\bvolume{32}(\bissue{6}),
\bfpage{3228}--\blpage{3250}
(\byear{2010})
\end{barticle}
\endbibitem

\bibitem[\protect\citeauthoryear{Shen and Wang}{2010}]{SW2010}
\begin{barticle}
\bauthor{\bsnm{Shen}, \binits{J.}},
\bauthor{\bsnm{Wang}, \binits{L.-L.}}:
\batitle{Sparse spectral approximations of high-dimensional problems based on
  hyperbolic cross}.
\bjtitle{SIAM Journal on Numerical Analysis}
\bvolume{48}(\bissue{3}),
\bfpage{1087}--\blpage{1109}
(\byear{2010})
\end{barticle}
\endbibitem

\bibitem[\protect\citeauthoryear{Harbrecht et~al.}{2008}]{Harbrecht2008NM}
\begin{barticle}
\bauthor{\bsnm{Harbrecht}, \binits{H.}},
\bauthor{\bsnm{Schneider}, \binits{R.}},
\bauthor{\bsnm{Schwab}, \binits{C.}}:
\batitle{Multilevel frames for sparse tensor product spaces}.
\bjtitle{Numerische Mathematik}
\bvolume{110}(\bissue{2}),
\bfpage{199}--\blpage{220}
(\byear{2008})
\end{barticle}
\endbibitem

\bibitem[\protect\citeauthoryear{Garcke et~al.}{2006}]{Garcke2006}
\begin{botherref}
\oauthor{\bsnm{Garcke}, \binits{J.}}, et al.:
Sparse grid tutorial.
Mathematical Sciences Institute, Australian National University, Canberra
  Australia
\textbf{7}
(2006)
\end{botherref}
\endbibitem

\bibitem[\protect\citeauthoryear{Zhang et~al.}{2013}]{zhang2013}
\begin{barticle}
\bauthor{\bsnm{Zhang}, \binits{G.}},
\bauthor{\bsnm{Gunzburger}, \binits{M.}},
\bauthor{\bsnm{Zhao}, \binits{W.}}:
\batitle{A sparse-grid method for multi-dimensional backward stochastic
  differential equations}.
\bjtitle{J. Comput. Math.}
\bvolume{31}(\bissue{3}),
\bfpage{221}--\blpage{248}
(\byear{2013})
\end{barticle}
\endbibitem

\bibitem[\protect\citeauthoryear{Zhang et~al.}{2016}]{Zhang2016NM}
\begin{barticle}
\bauthor{\bsnm{Zhang}, \binits{G.}},
\bauthor{\bsnm{Webster}, \binits{C.G.}},
\bauthor{\bsnm{Gunzburger}, \binits{M.}},
\bauthor{\bsnm{Burkardt}, \binits{J.}}:
\batitle{Hyperspherical sparse approximation techniques for high-dimensional
  discontinuity detection}.
\bjtitle{SIAM Review}
\bvolume{58}(\bissue{3}),
\bfpage{517}--\blpage{551}
(\byear{2016})
\end{barticle}
\endbibitem

\bibitem[\protect\citeauthoryear{Adcock et~al.}{2022}]{Adcock2022NM}
\begin{bbook}
\bauthor{\bsnm{Adcock}, \binits{B.}},
\bauthor{\bsnm{Brugiapaglia}, \binits{S.}},
\bauthor{\bsnm{Webster}, \binits{C.G.}}:
\bbtitle{Sparse Polynomial Approximation of High-dimensional Functions}
vol. \bseriesno{25},
(\byear{2022}).
\bcomment{SIAM}
\end{bbook}
\endbibitem

\bibitem[\protect\citeauthoryear{Barthelmann et~al.}{2000}]{Barthelmann2000}
\begin{barticle}
\bauthor{\bsnm{Barthelmann}, \binits{V.}},
\bauthor{\bsnm{Novak}, \binits{E.}},
\bauthor{\bsnm{Ritter}, \binits{K.}}:
\batitle{High dimensional polynomial interpolation on sparse grids}.
\bjtitle{Advances in Computational Mathematics}
\bvolume{12},
\bfpage{273}--\blpage{288}
(\byear{2000})
\end{barticle}
\endbibitem

\bibitem[\protect\citeauthoryear{Aronson}{1967}]{Aronson1967}
\begin{botherref}
\oauthor{\bsnm{Aronson}, \binits{D.G.}}:
Bounds for the fundamental solution of a parabolic equation
(1967)
\end{botherref}
\endbibitem

\bibitem[\protect\citeauthoryear{Sheng et~al.}{2024}]{SSX2025}
\begin{botherref}
\oauthor{\bsnm{Sheng}, \binits{C.}},
\oauthor{\bsnm{Su}, \binits{B.}},
\oauthor{\bsnm{Xu}, \binits{C.}}:
An implicit-explicit {M}onte {C}arlo method for semi-linear {PDE}s driven by
  $\alpha$-stable {L}{\'e}vy process and its error estimates.
Mathematics of Computation
(2024)
\end{botherref}
\endbibitem

\end{thebibliography}

\end{document}